%% file: ms.tex
\documentclass[11pt, a4paper, oneside,reqno]{amsart}

\input{AMIN_style}

\graphicspath{{Fig/}}

\usepackage[square,numbers]{natbib}

\newcommand{\bo}{T}

\newcommand{\sbo}{\widehat{T}}

\newcommand{\gr}{\nabla}
\newcommand{\pr}{_{\mathrm{prior}}}
\newcommand{\e}{\boldsymbol{1}}
\newcommand{\z}{\boldsymbol{0}}

\title[]{From Optimization to Control: \\ Quasi Policy Iteration}

\author[]{Mohamad Amin Sharifi Kolarijani$^{{\rm a}}$ and Peyman {Mohajerin Esfahani}$^{{\rm a,b}}$}
\thanks{The authors are with (a) Delft University of Technology, The Netherlands, and (b) University of Toronto, Canada.}
\thanks{This work was partially supported by the European Research Council (ERC) under the grant TRUST-949796, and the NSERC Discovery grant RGPIN-2025-06544.} 

\begin{document}

\begin{abstract}
Recent control algorithms for Markov decision processes (MDPs) have been designed using an implicit analogy with well-established optimization algorithms. 
In this paper, 
we adopt the quasi-Newton method (QNM) from convex optimization to introduce a novel control algorithm coined as quasi-policy iteration (QPI). 
In particular, QPI is based on a novel approximation of the ``Hessian'' matrix in the policy iteration algorithm, which exploits two linear structural constraints specific to MDPs and allows for the incorporation of prior information on the transition probability kernel. 
While the proposed algorithm has the same computational complexity as value iteration, it exhibits an empirical convergence behavior similar to that of QNM with a low sensitivity to the discount factor. 


\end{abstract}

\maketitle

%===============================================================================
\section{Introduction}
\label{sec:intro}
%===============================================================================
The problem of control, or the decision-making problem as it is also known within the operations research community, has been the subject of much research since the introduction of the Bellman principle of optimality in the late 1950s~\cite{bellman1957markovian}. 
In particular, the connection between control algorithms for Markov decision processes (MDPs) and optimization algorithms has been noticed since the late 1970s~\cite{puterman1979convergence}: 
Value iteration (VI)~\cite{bellman1957markovian} can be seen as an instance of gradient descent (GD) algorithm, and policy iteration (PI)~\cite{howard1960dynamic} is an instance of the Newton method (NM).

More recent works have used the relationship mentioned above to develop new control algorithms, with faster convergence and/or lower complexity, inspired by their counterparts for solving optimization problems~\cite{grand2021convex, vieillard2019connections}. 
In case of \emph{model-based} (a.k.a.~planning, with access to the model of the MDP) algorithms, 
the combination of the VI algorithm with Polyak momentum~\cite{polyak1964some}, Nesterov acceleration~\cite{nesterov1983method}, and Anderson acceleration~\cite{anderson1965iterative} have been explored in~\cite{goyal2019first} and~\cite{zhang2020globally}. 
More recently, Halpern anchoring acceleration~\cite{halpern1967fixed} has been used to introduce the Anchored VI algorithm~\cite{lee2024accelerating}, which in particular exhibits an~$\ord(1/k)$-rate for large values of discount factor and even for~$\gamma = 1$. 
Also of notice is the Generalized Second-Order VI algorithm~\cite{kamanchi2021generalized} which applies NM on a \emph{smoothed} version of the Bellman operator.
For \emph{model-free} (a.k.a.~learning, with access to samples from the MDP) algorithms, Speedy Q-Learning~\cite{ghavamzadeh2011speedy}, Momentum Q-Learning~\cite{weng2020momentum}, and Nesterov Stochastic Approximation~\cite{devraj2019matrix} are among the algorithms that use the idea of momentum for accelerating stochastic GD algorithm~\cite{yang2016unified,kidambi2018insufficiency,liu2018accelerating,allen2017katyusha} in order to achieve a better rate of convergence compared to standard Q-learning (QL). 
Moreover, the Zap Q-Learning algorithm~\cite{devraj2017zap} can be thought of as a second-order learning algorithm which was inspired by the stochastic Newton-Raphson (SNR) algorithm~\cite{ruppert1985newton}. 

We note that all of the aforementioned connections have been focused on \emph{finite} state-action MDPs. 
When it comes to infinite (continuous) state-action spaces, except in special cases such as linear–quadratic regulators (LQR), one needs to resort to finite-dimensional approximation techniques for computational purposes. 
This approximation may be at the modeling level by aggregation (discretization) of the state and action spaces, which readily falls into the finite state-action MDP setting~\cite{ref:Bert_disc_75,ref:Powell_07}. 
Alternatively, one may directly approximate the value function via finite parametrization and minimizing (a proxy of) the residual of its fixed-point characterization based on the Bellman principle of optimality~\cite{ref:Bert_abstract, szepesvari2022algorithms}. 
Examples of such include linear parameterization~\cite{ref:Bert_temporal, Tsitsiklis97}, or nonlinear parameterization with, for instance, neural network architectures~\cite{ref:Bert_neuro-dynamic,SCHMIDHUBER2015,sutton2018reinforcement} or max-plus approximation~\cite{Bach20, GONCALVES2021109623, ref:FDP_TAC, ref:FDP_NeurIPS, McEn06, amin2025fitted}. 
We also note that there is an alternative characterization of the original function as the solution to an infinite-dimensional linear program~\cite{ref:HL1}, paving the way for approximation techniques via finite tractable convex optimization~\cite{ref:VanRoy_MP,ref:HL2,ref:Moh_SIOPT}. 
With this view of the literature, it is worth noting that one can cast almost all of these approximation techniques as the solution to a finite-dimensional fixed-point or convex optimization problem. 

\noindent\textbf{Main contribution.} 
Motivated by the connection between control and optimization algorithms,
we adopt the quasi-Newton method (QNM) from convex optimization to introduce the quasi-policy iteration (QPI) algorithm with the following distinct features:
\begin{itemize}
        \item[1)] {\bf {Hessian approximation via structural information}}: 
        QPI is based on a novel approximation of the ``Hessian'' matrix in the PI algorithm by exploiting two linear structural constraints specific to MDPs and by allowing for the incorporation of prior information on the transition probability kernel of the MDP (Theorem~\ref{thm:QPI}). 
        In the special case of incorporating a uniform prior for the transition kernel, QPI can be viewed as a modification of the standard VI using two novel directions with adaptive step-sizes (Corollary~\ref{cor:QPI-A}). 
        \item[2)] \textbf{Convergence and sensitivity to discount factor}: 
        The per-iteration computational complexity of QPI is the same as VI, and its linear convergence can be guaranteed by safeguarding against standard VI~(Theorem~\ref{thm:QPI}) or backtracking~(Lemma~\ref{lem:backtrack}).  
        However, in our numerical simulations with random and structured MDPs, QPI exhibits an empirical behavior similar to QNMs with the convergence rate being less sensitive to the discount factor (Figure~\ref{fig:mb_instance}).
        \item[3)] \textbf{Local superlinear convergence}: We provide a modified implementation of QPI which also incorporates the secant-type constraints in approximation of the Hessian to guarantee the local superlinear convergence (Theorem~\ref{thm:QPI-B superlinear}).         
        \item[4)] \textbf{Extension to model-free control (a.k.a.~RL)}: We also introduce the quasi-policy learning (QPL) algorithm, the stochastic version of QPI, as a novel model-free algorithm with guaranteed convergence and the same per-iteration complexity as standard Q-learning (QL) algorithm (Theorem~\ref{thm:QPL}). 
\end{itemize}

The paper is organized as follows. 
In Section~\ref{sec:CORL}, we describe the optimal control problem of MDPs.  
In Section~\ref{sec:QPI}, we introduce and analyze the model-based QPI algorithm and its model-free extension, the QPL algorithm. 
All the technical proofs are provided in Section~\ref{app:proofs}. 
The performance of these algorithms is then compared with multiple control algorithms via extensive numerical experiments in Section~\ref{sec:experiments}. 
Section~\ref{sec:conclusion} concludes the paper by providing some final remarks. 

\noindent\textbf{Notations.} For a vector $v \in \R^n$, we use $v(i)$ and $[v](i)$ to denote its $i$-th element. 
Similarly, $M(i,j)$ and $[M](i,j)$ denote the element in row $i$ and column $j$ of the matrix $M \in \R^{m\times n}$. 
We use $\cdot\tr$ to denote the transpose of a vector/matrix. 
We use $\norm{\cdot}_2$ and $\norm{\cdot}_{\infty}$ to denote the 2-norm and $\infty$-norm of a vector, respectively. 
We use $\norm{\cdot}_2$ and $\norm{\cdot}_F$ for the induced 2-norm and the Frobenius norm of a matrix, respectively. 
Let $x\sim\PP$ be a random variable with distribution $\PP$. We particularly use $\hat{x}\sim\PP$ to denote \emph{a sample of the random variable $x$} drawn from the distribution $\PP$. 
We use $\e$ and $\z$ to denote the all-one and all-zero vectors, respectively. 
$I$ and $E = \e \e^{\top}$ denote the identity and all-one matrices, respectively. 
We denote the $i$-th unit vector by $e_i$, that is, the vector with its $i$-th element equal to $1$ and all other elements equal to $0$. 

%===============================================================================
\section{Optimal control of MDPs}
\label{sec:CORL}
%===============================================================================
A common formulation of the control problem relies on the concept of \emph{Markov decision processes} (MDPs). 
MDPs are a powerful modeling framework for stochastic environments that can be controlled to minimize some measure of cost. 
An MDP is a tuple $(\set{S},\set{A},\PP,c,\gamma)$, where $\set{S}$ and $\set{A}$ are the state space and action space, respectively. 
The transition kernel $\PP$ encapsulates the state dynamics: 
for each triplet $(s,a,s^+) \in \set{S}\times\set{A}\times\set{S}$, it gives the probability $ \PP(s^+|s,a)$ of the transition to state $s^+$ given that the system is in state $s$ and the chosen control is $a$. 
The cost function $c: \set{S}\times\set{A} \ra \R$, bounded from below, represents the cost $c(s,a)$ of taking the control action $a$ while the system is in state $s$. 
The discount factor $\gamma \in (0,1)$ can be seen as a trade-off parameter between short- and long-term costs. 
\emph{In this study, we consider tabular MDPs with a finite state-action space. 
In particular, we take $\set{S} = \{1,2,\ldots,n \}$ and $\set{A} = \{1,2,\ldots,m \}$.} 
This, in turn, allows us to treat functions $f:\set{S}\ra\R$ and $g:\set{S}\times\set{A}\ra\R$ as vectors 
$f\in \R^{|\set{S}|} = \R^{n}$ and $g\in \R^{|\set{S}\times\set{A}|} = \R^{nm}$ -- in the latter case, we are considering a proper 1-to-1 mapping $\set{S}\times\set{A} \ra \{1,2,\ldots,nm \}$. 

Let us now fix a control policy $\pi:\set{S}\ra\set{A}$, i.e., a mapping from states to actions. 
The stage cost of the policy $\pi$ is denoted by $c^{\pi} \in \R^{|\set{S}|} = \R^{n}$ with elements $c^{\pi}(s) = c\big(s,\pi(s)\big)$ for $s \in \set{S}$. 
The transition (probability) kernel of the resulting Markov chain under the policy $\pi$ is denoted by $\PP^{\pi}$, where $\PP^{\pi} (s^+|s) = \PP\big(s^+|s,\pi(s)\big)$ for $s,s^+\in \set{S}$. 
We also define the matrix $P^\pi \in \R^{|\set{S}|\times|\set{S}|} = \R^{n\times n}$, with elements $P^\pi(s,s^+) \Let \PP^{\pi}(s^+ | s)$ for $s,s^+ \in \set{S}$, to be the corresponding transition (probability) \emph{matrix}. 
The value of a policy is the expected, discounted, accumulative cost of following this policy over an infinite-horizon trajectory:  
For the policy $\pi$, we define the value function $v^{\pi} \in \R^{|\set{S}|} = \R^n$ with elements 
$$v^{\pi}(s) \Let \mathbb{E}_{s_{t+1}\sim\mathbb{P}^{\pi}(\cdot|s_t)}\big[\ssum_{t=0}^\infty\gamma^t c^{\pi}(s_t) \mid s_0 = s\big],$$ 
and the Q-function $q^{\pi} \in \R^{|\set{S}\times\set{A}|} = \R^{nm}$ with elements 
$$q^{\pi}(s,a) \Let c(s,a) + \gamma \mathbb{E}_{s^+\sim\mathbb{P}(\cdot|s,a)} [v^{\pi}(s^+)],$$ 
so that we also have $v^\pi(s) = q^\pi\big((s,\pi(s)\big)$ for each $s\in\set{S}$. 
Given a value function $v$, let us also define $\pi_v:\set{S}\ra\set{A}$ by $$\pi_v(s) \in \argmin_{a\in\set{A}}\left\{ c(s,a) + \gamma \EE_{s^+\sim\PP(\cdot|s,a)}\left[v(s^+)\right]\right\},$$ 
to be the \emph{greedy policy w.r.t.~$v$}. 
Similarly, for a Q-function~$q$, define $\pi_q:\set{S}\ra\set{A}$ by 
$$\pi_q(s) \in \argmin_{a\in\set{A}} q(s,a),$$ 
to be the \emph{greedy policy w.r.t.~$q$}. 
The problem of interest is to control the MDP optimally, that is, to find the optimal policy $\pi^*$ with the optimal value/Q-function 
\begin{equation}\label{eq:optimal value func}
    v\opt = \min_{\pi} v^{\pi}, \quad q\opt = \min_{\pi} q^{\pi},
\end{equation}
so that the expected, discounted, infinite-horizon cost is minimized. 
Let us also note that the optimal policy, i.e., the minimizer of the preceding optimization problems, is the greedy policy w.r.t.~$v\opt$ and $q\opt$, that is, $\pi\opt = \pi_{v\opt} = \pi_{q\opt}$. 

Interestingly, the optimal value/Q-function introduced in~\eqref{eq:optimal value func} can be equivalently characterized as the fixed-point of the corresponding Bellman operators. 
This fixed-point characterization is the basis for the class of value iteration algorithms. 
To be precise, we have $ v\opt = \bo(v\opt)$, i.e., 
\begin{equation}\label{eq:fp_v}
    v\opt(s) = [\bo (v\opt)](s), \quad \forall s \in \set{S},
\end{equation}
where $\bo:\R^{|\set{S}|}\ra\R^{|\set{S}|}$ is the \emph{Bellman operator} given by 
\begin{equation} \label{eq:T}
    [\bo (v)](s) \Let \min_{a \in\set{A}} \left\{ c(s,a) + \gamma \EE_{s^+\sim\PP(\cdot|s,a)} \left[v(s^+)\right] \right\}. 
\end{equation}
Similarly, we have $q\opt  =  \EE[\sbo (q\opt)]$, i.e., 
\begin{equation}\label{eq:fp_q}
            q\opt(s,a) = \EE\left[[\sbo (q\opt)](s,a) \right] ,\quad \forall (s,a)\in\set{S}\times\set{A}
\end{equation}
where $\sbo:\R^{|\set{S}\times\set{A}|}\ra\R^{|\set{S}\times\set{A}|}$ is the \emph{sampled} Bellman operator given by\footnote{Strictly speaking, the provided sampled Bellman operator is the empirical version of the Bellman operator for the Q-function, given by $[\bar{\bo} (q)](s,a) \Let c(s,a) + \gamma \EE_{s^+\sim\PP(\cdot|s,a)} \left[\min_{a^+ \in\set{A}} q(s^+,a^+)\right]$ for $(s,a)\in\set{S}\times\set{A}$.}
\begin{equation}\label{eq:rT}
    [\sbo (q)](s,a) \Let c(s,a) + \gamma \min_{a^+ \in \set{A}} q(\hat{s}^+,a^+), 
\end{equation}
with $\hat{s}^+ \sim \PP(\cdot|s,a)$ being a \emph{sample} of the next state drawn from the distribution $\PP(\cdot|s,a)$ for the pair $(s,a)$. 
We note that the expectation in~\eqref{eq:fp_q} is w.r.t.~$\hat{s}^+$. 

%===============================================================================
\section{Quasi-Policy Iteration (QPI)}\label{sec:QPI}
%===============================================================================

In this section, we develop and analyze the quasi-policy iteration (QPI) algorithm for model-based control problems and its extension, the quasi-policy learning (QPL) algorithm, for model-free control problems.\footnote{In this paper, the terminologies of ``model-free'' and ``model-based'' indicate the {\em available information (oracle)}, i.e., whether we have access to the model or only the system trajectory (samples). We note that this is different from the common terminologies in the RL literature where these terms refer to the {\em solution approach}, i.e., whether we identify the model along the way (model-based RL) or directly solve the Bellman equation to find the value function (model-free RL).} 
These algorithms, as the name suggests, are inspired by quasi-Newton methods (QNMs) from convex optimization. 
For this reason, we begin with a brief overview of the connection between optimization and control algorithms, as well as QNMs, to motivate the algorithms developed in this study. 

We note that for tabular MDPs with $\set{S} = \{1,\ldots,n\}$ and $\set{A} = \{1,\ldots,m\}$, we have $v\in \R^{n}$ and $q\in \R^{nm}$ for the value function and the Q-function, respectively. 
Correspondingly, we have $T : \R^n\ra\R^n$ with
\begin{align*}
    \bo (v) = \ssum_{s \in\set{S}} \ [\bo (v)](s) \cdot e_s,
\end{align*} 
where $e_s \in \R^{n}$ is the unit vector for the state $s \in \set{S}$, 
and also $\sbo : \R^{nm} \ra \R^{nm}$ with
\begin{align*}
    \sbo (q) = \ssum_{(s,a) \in\set{S}\times\set{A}} \ [\sbo (q)](s,a) \cdot e_{(s,a)},
\end{align*}
where $e_{(s,a)} \in \R^{nm}$ is the unit vector for the state-action pair~$(s,a) \in \set{S}\times\set{A}$. 

\subsection{Optimization vs. Control}

Two classical algorithms for solving the optimal control problem of MDPs are value iteration (VI) and policy iteration (PI), which are instances of gradient descent (GD) and Newton method (NM), respectively, as we mentioned above. 
To be precise, consider the unconstrained minimization problem 
\[ 
\min_{x \in \R^{\ell}} f(x),
\]
where $f:\R^{\ell}\ra\R$ is twice continuously differentiable and strongly convex, and let $g(x) = \gr f(x)$ and $H(x) = \gr^2f(x)$ be the gradient and Hessian of $f$ evaluated at $x$. 
Then, by defining $g(v) \Let v - \bo (v)$, one can see that the VI update rule~\cite{bellman1957markovian}
\[
v_{k+1} = \bo (v_k) = v_k - g(v_k),
\]
resembles the GD update rule 
$$x_{k+1} = x_k - \alpha_k g(x_k),$$ 
where $\alpha_k$ is the step-size. 
Similarly, by defining $H(v) \Let I - \gamma P^{\pi_{v}}$, the PI update rule~\cite[Prop.~6.5.1]{puterman2014markov}
\[
v_{k+1} = \big(I - \gamma P^{\pi_{v_k}}\big)^{-1} c^{\pi_{v_k}} = v_k - [H(v_k)]^{-1} g(v_k),
\]
resembles the NM update rule 
$$x_{k+1} = x_k - \alpha_k [H(x_k)]^{-1}g(x_k),$$ 
where $\alpha_k$ is again the step-size. 

The VI algorithm converges linearly with rate $\gamma$ owing to the fact that the operator~$\bo$ is a~$\gamma$-contraction in the~$\infty$-norm~\cite[Prop.~6.2.4]{puterman2014markov}. 
The PI algorithm, on the other hand,  outputs the optimal policy in a finite number of iterations~\cite[Thm.~6.4.2]{puterman2014markov}.
Moreover, the algorithm has a local \emph{quadratic} rate of convergence when initiated in a small enough neighborhood around the optimal solution~\cite{bertsekas2022lessons, gargiani2022dynamic}. 
The faster convergence of PI compared to VI, however, comes with a higher per-iteration computational complexity: The per-iteration complexities of VI and PI are~$\ord(n^2m)$ and~$\ord(n^2m+n^3)$, respectively. 
The extra~$\ord(n^3)$ complexity is due to the policy evaluation step, i.e., solving a linear system of equations; see also the matrix inversion in the characterization above. 
Not surprisingly, we see the same convergence-complexity trade-off between GD and NM. 
While NM has a better convergence rate compared to GD, it suffers from a higher per-iteration computational cost. 
GD, with a proper choice of step-size, converges \emph{linearly}~\cite[Thm. 3.12]{bubeck2015convex} with $O(\ell)$ per-iteration complexity (disregarding the complexity of gradient oracle). 
On the other hand, NM, with a proper choice of step-size, has a local \emph{quadratic} convergence rate~\cite[Thm.~5.3]{bubeck2015convex} with $O(\ell^3)$ per-iteration complexity, assuming direct inversion (and disregarding the complexity of gradient and Hessian oracles). 

QNMs are a class of methods that allow for a trade-off between computational complexity and (local) convergence rate. 
To do so, these methods use a Newton-type update rule 
$$x_{k+1} = x_k - \alpha_k \wt{H}_k^{-1}\gr f(x_k),$$ 
where $\wt{H}_k$ is an \emph{approximation} of the true Hessian $\gr^2f(x_k)$ at iteration $k$. 
Different QNMs use different approximations of the Hessian. 
A generic approximation scheme in QNMs is 
\begin{equation}\label{eq:QNM}
\begin{array}{rl}
    \wt{H}_k = \argmin\limits_{H \in \R^{\ell \times \ell}} & \norm{H-H\pr}^2_F \\
               \text{s.t.} & H r_i = b_i,\ i =1,\ldots,j,
\end{array}
\end{equation}
which minimizes the distance (in Frobenius norm) to a given prior~$H\pr$  subject to $j$ ($\geq 1$) linear constraints specified by $r_i,b_i \in \R^{\ell}$. 
This leads to the approximation~$\wt{H}_k$ being a rank-$j$ update of the prior~$H\pr$, i.e.,
\begin{equation}\label{eq:QNM solution}
    \wt{H}_k = H\pr + (B - H\pr R)(R\tr R)^{-1} R\tr,
\end{equation}
where $R = (r_1, \ldots, r_j),\ B = (b_1, \ldots, b_j) \in \R^{\ell \times j}$. Hence, $\wt{H}_k^{-1}$ can be easily computed based on $H\pr^{-1}$ using the Woodbury formula.
Different choices of the prior and the linear constraints in the generic approximation scheme above lead to different QNMs. 
For example, by choosing the so-called \emph{secant conditions} with $r_i = x_{k-i+1} - x_{k-i}$ and $b_i = \gr f(x_{k-i+1}) - \gr f(x_{k-i})$ as linear constraints and $H\pr=I$ as the prior, we derive Anderson mixing with memory~$j$~\cite{anderson1965iterative}, while by using a single secant condition with $j=1$ and choosing $H\pr=\wt{H}_{k-1}$ as the prior, we derive QNM with Broyden approximation~\cite{broyden1965class}.

In what follows, we use a similar idea and propose the \emph{quasi-policy iteration (QPI)} algorithm by incorporating a computationally efficient approximation of the ``Hessian'' $H = I - \gamma P$ in the PI algorithm. 
We note that the authors in~\cite{geist2018anderson, zhang2020globally, sun2021damped} also propose the combination of Anderson mixing with optimal control algorithms. 
However, the QPI algorithm is fundamentally different in the sense that it approximates the transition matrix $P$ using a different set of constraints that are specific to the optimal control algorithms.

\subsection{QPI Algorithm} 
For $k\in \{0,1,2,\ldots\}$, let 
$$
c_k\Let c^{\pi_{v_k}},\quad P_k \Let P^{\pi_{v_k}},\quad \bo_k \Let \bo(v_k).
$$
(Recall that $c^{\pi_{v_k}}$ and $P^{\pi_{v_k}}$ are the stage cost and the state transition matrix of the greedy policy $\pi_{v_k}$ w.r.t.~$v_k$, respectively.) 
Recall the PI update rule
\begin{align*}
    v_{k+1} &= (I - \gamma P_k)^{-1} c_k 
    =  v_k - (I - \gamma P_k)^{-1} (  v_k - \bo_k ).
\end{align*}
Inspired by the QNM approximation scheme~\eqref{eq:QNM}, we propose the generic QPI update rule
\begin{equation}\label{eq:QPI update general}
    v_{k+1} =  v_k - (I - \gamma \wt{P}_k)^{-1} (  v_k - \bo_k ),
\end{equation}
where
\begin{equation}\label{eq:QPI approx general}
\begin{array}{rl}
    \wt{P}_k = \argmin\limits_{P \in \R^{n \times n}} & \norm{P-P\pr}^2_F \\
               \text{s.t.} & P r_i = b_i,\ i =1,\ldots,j.
\end{array}
\end{equation}
Observe that instead of approximating the complete Hessian $H_k\Let (I - \gamma P_k)^{-1}$ similar to standard QNMs, we are only approximating $P_k$. 
This choice particularly allows us to exploit the problem structure in order to form novel constraints and prior as we discuss next.

Regarding the constraints, the problem structure gives us two linear equality constraints: 
First, $P_k$ is a row stochastic matrix, i.e., we have the \emph{global} structural constraint
\begin{equation}\label{eq:normalize_cond}
    P_k  \e =  \e, 
\end{equation}
and hence we can set $r_1 = b_1 = \e$. 
Second, we can use the fact that the Bellman operator~$\bo$ is piece-wise affine. 
In particular, from the definition~\eqref{eq:T} of the Bellman operator, it follows that $\bo(v) = c^{\pi_v} + \gamma P^{\pi_v} v$. Thus, we have the \emph{local} structural constraint
\begin{equation}\label{eq:affine_cond}
    \bo_k =  c_k + \gamma P_k v_k \; \Rightarrow \; P_k v_k = {\gamma}^{-1}  (  \bo_k -  c_k ),
\end{equation} 
and we can set $r_2 = v_k$ and $b_2 = {\gamma}^{-1}  (  \bo_k -  c_k )$. 
Note that, unlike the standard secant conditions in QNMs, the constraints~\eqref{eq:normalize_cond} and~\eqref{eq:affine_cond} hold \emph{exactly}. 
Incorporating these constraints, we propose the approximation
\begin{equation}\label{eq:QPI approx}
\begin{array}{rl}
    \wt{P}_k = \argmin\limits_{P \in \R^{n \times n}} & \Vert P - P\pr \Vert^2_F \\
               \text{s.t.} & P \e = \e,\
               P v_k = {\gamma}^{-1} (  \bo_k -  c_k ).
\end{array}
\end{equation}

The update rule~\eqref{eq:QPI update general} using the approximation~\eqref{eq:QPI approx} is, however, not necessarily a contraction. 
The same problem also arises in similar algorithms such as Anderson accelerated VI~\cite{zhang2020globally} and Nesterov accelerated VI~\cite{goyal2019first}. 
Here, we follow the standard solution for this problem, that is, safeguarding the QPI update against the standard VI update based on the Bellman error 
$$\theta_k \Let \norm{v_{k} - \bo_{k} }_{\infty}.$$ 
To be precise, at each iteration $k=0, 1,\ldots$, we consider the \emph{safeguarded} QPI update rule as follows
\begin{equation}\label{eq:QPI safegaurd}
\begin{array}{ll}
    \text{(QPI)} & \text{compute $v_{k+1}$ according to~\eqref{eq:QPI update general},~\eqref{eq:QPI approx}}; \\
    \text{(Safeguard)} & \text{{\bf if} $\theta_{k+1} > \gamma^{k+1} \theta_0$, {\bf then} $v_{k+1} = \bo_k$.} 
\end{array}
\end{equation}

The following theorem summarizes the discussion above by providing the QPI update rule explicitly (see Section~\ref{appendix:proof thm} for the proof). 

\begin{Thm}[QPI] \label{thm:QPI}
Consider the update rule~\eqref{eq:QPI update general} using the approximation~\eqref{eq:QPI approx} where $P\pr \e = \e$ and let $G\pr = (I-\gamma P\pr)^{-1}$. 
We have
\begin{subequations} \label{eq:QPI update gen}
\begin{equation}\label{eq:QPI update}
     v_{k+1} = v_k -  \wt{G}_{k} (v_k - \bo_k),
\end{equation}
where
\begin{equation}\label{eq:QPI-def0}
    \begin{array}{l}
        w_k = \bo_k -c_k - \gamma P\pr v_k, \ \check{w}_k = G\pr w_k \in \R^{n}, \\[1ex]
        u_k = v_k - \frac{\e\tr v_k}{n} \e,\ \check{u}_k = G\pr\tr u_k \in \R^{n}, \\[1ex]
        \tau_k = \left\{ \begin{array}{ll}
        0 &  \text{if}\  u_k\tr v_k = 0, \\
        (u_k\tr v_k)^{-1}  & \text{otherwise}
        \end{array}\right. \in \R,\\[2ex]
        \eta_k = \left\{ \begin{array}{ll}
        0 & \text{if}\  u_k\tr v_k = 0, \\
        \big(u_k\tr (v_k - \check{w}_k) \big)^{-1}  & \text{otherwise}
        \end{array}\right. \in \R,\\[2ex]
        \wt{P}_k = P\pr + \gamma^{-1} \tau_k w_k u_k\tr, \\
        \wt{G}_k  =  G\pr + \eta_k \check{w}_k \check{u}_k\tr.
    \end{array}
\end{equation}
\end{subequations}
Moreover, the iterates $v_k$ of the QPI update rule~\eqref{eq:QPI update gen} with the safeguarding~\eqref{eq:QPI safegaurd} converge to $v\opt$ linearly with rate $\gamma$ and a per-iteration time complexity of $\ord(n^2m)$.
\end{Thm} 

Observe that the safeguarded QPI update rule has the same per-iteration complexity as VI. 
Moreover, the convergence of QPI is ensured via safeguarding against VI, which leads to the same theoretically guaranteed linear convergence with rate~$\gamma$ as for VI.
However, as we will show in the numerical examples below, we observe an empirically faster convergence for QPI with its rate showing less sensitivity to $\gamma$ similar to QNMs. 
We also note that there is also a one-time computational cost of $O(n^3)$ in the QPI update rule ~\eqref{eq:QPI update gen} for computing $G\pr$ (assuming direct inversion and if $G\pr$ is not available in closed form). 

\subsubsection{The prior} 
Next to be addressed is the choice of~$P\pr$. 
First, note that for a fixed prior in all iterations, computation of $\tau_k$ and $\wt{P}_k$ in~\eqref{eq:QPI-def0} is not needed since the update~\eqref{eq:QPI update} only requires $\wt{G}_k$. 
The first choice for such a fixed prior is to exploit the available knowledge on the structure of the MDP. 
For instance, one can set $P\pr = P^{\mu}$ with $\mu$ being the stochastic policy choosing actions uniformly at random (so that $P\pr$ is the average over actions of and has the same sparsity pattern as the true transition kernel of the MDP).  
A computationally advantageous choice for the prior is the uniform distribution~$ P\pr = \frac{1}{n} E = \frac{1}{n}\e \e\tr$ for which the update rule can be simplified significantly (see Section~\ref{appendix:proof cor} for the proof):

\begin{Cor}[Uniform prior] \label{cor:QPI-A} 
The QPI update rule~\eqref{eq:QPI update gen} with the uniform prior~$ P\pr = \frac{1}{n} E$,  equivalently reads as 
\begin{subequations} \label{eq:QPI-A update gen}
\begin{equation}\label{eq:QPI-A update}
     v_{k+1} = (1-\delta_k) \bo_k + \delta_k c_k + \lambda_k \e,
\end{equation}
where the scalar coefficients are given by  
\begin{equation}\label{eq:QPI-A-def0}
    \begin{array}{l}
         z_k = c_k - \frac{\e\tr c_k}{n} \e  \in \R^n,\\
         g_k = v_k - \bo_k \in \R^n,\quad 
         y_k = g_k  - \frac{\e\tr g_k }{n} \e \in \R^n, \\[1ex]
         \delta_k = \left\{ \begin{array}{ll}
             0 & \text{if}\  v_k\tr (y_k+z_k) = 0, \\
            \frac{v_k\tr y_k}{v_k\tr (y_k+z_k)}  & \text{otherwise}
         \end{array} 
         \right. \in \R,\\[3ex]
         \lambda_k = \frac{\gamma}{n(1-\gamma)} \e\tr \big((\delta_k-1) g_k + \delta_k c_k\big) \in \R. %\\[2ex]
    \end{array}
\end{equation}
\end{subequations}
\end{Cor}

Observe that the QPI update rule~\eqref{eq:QPI-A update} is a modification of the standard VI update rule $v_{k+1} = \bo (v_k)$ using two new vectors, namely, $c_k - \bo (v_k)$ and the all-one vector~$\e$, with adaptive coefficients $\delta_k$ and $\lambda_k$, respectively. 

Another interesting choice for the prior in~\eqref{eq:QPI approx} is $P\pr = \wt{P}_{k-1}$, i.e., the previous approximation. 
This leads to a recursive scheme for approximating the transition matrix similar to QNM with Broyden approximation~\cite{broyden1965class}. 
This can be achieved by choosing an initialization~$\wt{P}_{-1}$ such that $\wt{P}_{-1} \e = \e$ and defining~$\wt{G}_{-1} \Let (1-\gamma \wt{P}_{-1})^{-1}$ (e.g., $\wt{P}_{-1} = \frac{1}{n}E$ and $\wt{G}_{-1} = I + \frac{\gamma}{n(1-\gamma)}E$). 
We note that for this recursive scheme, one might also need to introduce some form of projection to avoid unbounded growth of $\wt{P}_{k}$ and $\wt{G}_{k}$. 
In particular, using the fact that $\wt{P}_{k}$ is in principle an approximation of a row stochastic matrix, one can use a simple normalization to ensure $\|\wt{P}_{k}\|_\infty \leq 1$ and $\|\wt{G}_{k}\|_\infty \leq (1-\gamma)^{-1}$ for all $k$.

\subsubsection{Implementation via backtracking}\label{sec:QPI backtrack}
The QPI update rule~\eqref{eq:QPI update} can be alternatively implemented without safeguarding. The proposed alternative again exploits the fact that the VI update leads to contraction in the Bellman error and combines that with the \emph{backtracking} approach from optimization. 
To be precise, let us rewrite the QPI update rule~\eqref{eq:QPI update} as
\begin{align*}
     v_{k+1} &= v_k -  (\wt{G}_{k}+ I - I) (v_k - \bo_k) 
     = \bo_k - (\wt{G}_{k} - I) (v_k - \bo_k),
\end{align*}
and introduce a step-size $\alpha_k$ as follows 
\begin{equation}\label{eq:QPI backtrack}
     v_{k+1} = \bo_k - \alpha_k (\wt{G}_{k} - I) (v_k - \bo_k).
\end{equation}
Then, at each iteration $k$, we introduce an inner iteration that finds the step-size~$\alpha_k$ via backtracking such that $\theta_{k+1} \leq \gamma' \theta_{k}$ for a chosen $\gamma' \in (\gamma, 1)$. 
For instance, at each iteration~$k$,  we can set
\begin{equation}\label{eq:QPI backtrack 1}
\begin{array}{ll}
    \text{(0)} & \alpha_k = 1; \\
    \text{(1)} & \text{compute $v_{k+1}$ according to~\eqref{eq:QPI backtrack},~\eqref{eq:QPI-def0}}; \\
    \text{(2)} & \text{{\bf if} $\theta_{k+1} > \gamma' \theta_k$, {\bf then} $\alpha_k \gets \frac{1}{2} \alpha_k$, go to (1).} 
\end{array}
\end{equation}
We note that the preceding backtracking scheme is guaranteed to terminate if $v_k \neq v\opt$ (i.e, $\theta_k \neq 0$). 
Indeed, we have (see Section~\ref{appendix:proof lem backtrack} for the proof):
\begin{Lem}[Backtracking]\label{lem:backtrack}
    Let $\theta_v \Let \norm{v-\bo(v)}_{\infty}$ for $v\in\R^n$. 
    Fix $v\in\R^n$ such that $\theta_v \neq 0$. 
    Then, for all $w\in\R^n$ with $\norm{w-\bo(v)}_{\infty} \leq \frac{\gamma'-\gamma}{1+\gamma}\theta_v$, we have $\theta_w \leq \gamma' \theta_v$. 
\end{Lem}
In particular, the preceding lemma implies that the backtracking scheme~\eqref{eq:QPI backtrack 1} terminates in at most $\log_2(\frac{1+\gamma}{\gamma'-\gamma}\Vert\wt{G}_{k} - I\Vert_{\infty})$ steps at iteration~$k$. 
We also note that iterates of the proposed QPI update rule with backtracking converge to the optimal value function linearly with rate~$\gamma'$. 

\subsubsection{Modified implementation with superlinear convergence}\label{sec:superlinear} 
For the proposed QPI scheme to achieve the classic local \emph{superlinear} convergence of the class of QNMs, we need to modify the approximation~\eqref{eq:QPI approx} to also include secant-type conditions. 
To this end, for $k\geq 1$, we consider a modified approximations as follows

\begin{subequations} \label{eq:QPI-B approx}
\text{{\bf if} $\pi_{v_k} = \pi_{v_{k-1}}$}
\begin{equation}\label{eq:QPI-B selection 1}
\begin{array}{rl}
    \wt{P}_k = \argmin\limits_{P \in \R^{n \times n}} & \Vert P - \wt{P}_{k-1} \Vert^2_F \\
    \text{s.t.} & P \e = \e,\
     P v_k = {\gamma}^{-1} (  \bo_k -  c_k ), \
     \gamma P (v_k-v_{k-1}) = \bo_k-\bo_{k-1};
\end{array} 
\end{equation}

\text{{\bf else}}
\begin{equation}\label{eq:QPI-B selection 2}
\begin{array}{rl}
    \wt{P}_k = \argmin\limits_{P \in \R^{n \times n}} & \Vert P - \wt{P}_{k-1} \Vert^2_F \\
    \text{s.t.} & P \e = \e, \
     \gamma P (v_k-v_{k-1}) = \bo_k-\bo_{k-1};
\end{array}
\end{equation}
with the corresponding gain matrix
\begin{equation}\label{eq:QPI-B gain}
\wt{G}_k = \big(I -\gamma \wt{P}_k\big)^{-1}.
\end{equation}
\end{subequations}
Above, the only difference between \eqref{eq:QPI-B selection 1} and \eqref{eq:QPI-B selection 2} is the omission of the \emph{local structural} constraint in \eqref{eq:QPI-B selection 2}. 
We note that the condition $\pi_{v_k} = \pi_{v_{k-1}}$ implies that $v_k$ and $v_{k-1}$ have the same greedy policy and hence lie on the same affine region of the map $T$ (recall that $T$ is piece-wise affine). 
Therefore, the selection rule is to avoid ``constraint collision'' between the \emph{local structural} and the \emph{secant} constraints when $\pi_{v_k} \neq \pi_{v_{k-1}}$.
Also, observe that in both approximations the prior is the previous approximation~$\wt{P}_{k-1}$. 
Most importantly, note that the solution~$\wt{P}_k$ is a low-rank (rank-three, at most) update of~$\wt{P}_{k-1}$; c.f.~the optimization~\eqref{eq:QNM} and its solution~\eqref{eq:QNM solution}. 
This means that $\wt{G}_k$ can also be computed with a low cost using the Woodbury formula and $\wt{G}_{k-1} = (I -\gamma \wt{P}_{k-1})^{-1}$. 
Our next result concerns the convergence of the proposed modified implementation of the QPI algorithm (see Section~\ref{appendix:proof thm QPI-B superlinear} for the proof):

\begin{Thm}[Superlinear convergence] \label{thm:QPI-B superlinear}
Consider the update rule~\eqref{eq:QPI backtrack} with the gain matrix~\eqref{eq:QPI-B approx} and the backtracking scheme~\eqref{eq:QPI backtrack 1}. Then, the iterates $v_k$ converge to $v\opt$ with a global linear rate and a local superlinear rate. 
\end{Thm}

\subsubsection{Other constraints} 
We finish this section with the following remark.  
One can also add extra constraints to the minimization problem~\eqref{eq:QPI approx} to impose a particular structure on the approximate transition matrix $\wt{P}_k$. 
For instance, a natural constraint is to require this matrix to be entry-wise non-negative so that $\wt{P}_k$ is indeed a probability transition matrix; or, one can impose a sparsity pattern on $\wt{P}_k$ using the existing knowledge on the structure of the MDP. 
However, incorporating such information may lead to the problem~\eqref{eq:QPI approx} not having a closed-form and/or low-rank solution, and hence undermining the computational efficiency of the proposed algorithm. 
In this regard, we note that the problem~\eqref{eq:QPI approx} has a closed-form solution which is a rank-one update of the prior; see $\wt{P}_{k}$ in~\eqref{eq:QPI-def0}.

\subsection{Extension to model-free control: QPL algorithm} 
We now introduce the quasi-policy learning (QPL) algorithm as the extension of QPI for model-free control problems with access to samples through a generative model. 
For simplicity, we limit the following discussion to the extension of the QPI algorithm~\eqref{eq:QPI-A update gen} with a uniform prior. 
However, we note that the extension can be similarly applied to the QPI algorithm~\eqref{eq:QPI update gen} with the generic prior. 

The basic idea is to implement the stochastic version of the QPI update rule \emph{for the Q-function} using the samples. 
In particular, similar to the approximation~\eqref{eq:QPI approx}, we use an approximation of the state-action transition matrix under the greedy policy w.r.t.~the Q-function~$q_k$ at each iteration~$k$. 
However, the \emph{local} structural constraint is in this case formed based on the \emph{sampled} Bellman operator $\sbo (\cdot)$, evaluated at the sampled next states $\hat{s}_k^+$ at iteration $k$, as a surrogate for the Bellman operator $\bo (\cdot)$. 
To be precise, let 
$$\sbo_k \Let \sbo (q_k),$$ 
at each iteration $k$. 
Also, let $c \in \R^{nm}$ be the vector of stage cost (with the same state-action ordering as the Q-function $q_k \in\R^{nm}$). 
We note that since the proposed QPL algorithm is synchronous with one sample for each state-action pair in each iteration, we have access to the complete stage cost~$c$ after the first iteration and can treat it as an input to the algorithm. 
The approximate state-action transition matrix~$\wt{P}_k$ at each iteration~$k$ is then formed as follows
\begin{equation}\label{eq:QPL-approx}
\begin{array}{rl}
    \wt{P}_k = &\argmin\limits_{P \in \R^{nm \times nm}}  \Vert P - P\pr \Vert^2_F \\
               & \ \ \text{s.t.} \ \   P \e = \e,\ P q_k = {\gamma}^{-1} \big(  \sbo_k -  c \big).
\end{array}
\end{equation}
The minimization problem above also has a closed-form solution as a rank-one update of the prior, which allows us to compute $\wt{G}_k = (I-\gamma \wt{P}_k)^{-1}$ efficiently using  Woodbury formula. 
Having $\wt{G}_k$ at hand, we can consider the update rule
\begin{equation}\label{eq:QPL simple}
q_{k+1} = q_k - \alpha_k \wt{G}_k (q_k-\sbo_k),
\end{equation}
where $\alpha_k$ is the diminishing learning rate of the algorithm, e.g., $\alpha_k = 1/(k+1)$.
However, similar to the model-based case, the update rule~\eqref{eq:QPL simple} is not convergent. 
To address this issue, we consider a convex combination of the update rule~\eqref{eq:QPL simple} with the ``standard'' synchronous Q-learning (QL) update rule\footnote{This is the so-called \emph{synchronous} update of the Q-function in \emph{all} state-action pairs in each iteration, corresponding to the \emph{parallel sampling model} introduced by~\cite{kearns1998finite}.}~\cite{watkins1992q,kearns1998finite}
\begin{equation}\label{eq:QL}
q_{k+1} = q_k - \alpha_k (q_k-\sbo_k),
\end{equation}
with a diminishing effect. That is, we consider the update rule
\begin{align*}
q_{k+1} &= q_k - \alpha_k \bigg( (1-\beta_k) (q_k-\sbo_k) + \beta_k  \wt{G}_k (q_k-\sbo_k) \bigg)\\
&= q_k + \alpha_k (\sbo_k-q_k) + \alpha_k \beta_k  (\wt{G}_k-I) (\sbo_k-q_k),
\end{align*}
where $\beta_k$ is a diminishing coefficient, e.g., $\beta_k = 1/(k+1)$. 
Moreover, to assure convergence, we bound the effect of the extra term $p_k = (\wt{G}_k - I) (\sbo_k-q_k)$ by passing it through the operator~$\Pi_M: \R^{nm}\ra\R^{nm}$, with $M = \frac{2\gamma}{(1-\gamma)^2}\norm{c}_\infty$, given by
\begin{equation*}
  \Pi_M(p) \Let \frac{\min\{M, \norm{p}_\infty\}}{ \norm{p}_\infty} p.
\end{equation*} 
We note that the bound $M$ is chosen using the fact that for any row stochastic matrix $P$ and $G = (1-\gamma P)^{-1}$, we have 
\[
\norm{G - I}_\infty = \norm{\ssum_{i=1}^\infty (\gamma P)^{i}}_\infty \leq \ssum_{i=1}^\infty \gamma^i \norm{P^{i}}_\infty \leq \frac{\gamma}{1-\gamma},
\]
and the generic bound $\norm{q^{\pi}}_\infty \leq \norm{c}_\infty/(1-\gamma)$ for the value $q^{\pi}$ of any policy~$\pi$. 
The update rule of the model-free QPL algorithm in its generic form is hence
\begin{equation}\label{eq:QPL generic}
\begin{array}{l}
     p_k = (\wt{G}_k - I) (\sbo_k-q_k),  \\[1ex] 
     q_{k+1} = q_k + \alpha_k (\sbo_k-q_k) + \alpha_k \beta_k \Pi_M(p_k). 
\end{array}
\end{equation}
In particular, by using the uniform prior~$P\pr = \frac{1}{nm} E$, the QPL update rule reduces to
\begin{subequations}\label{eq:QPL update gen}
\begin{equation}\label{eq:QPL update}
    \begin{array}{l}
         p_k = \delta_k (c-\sbo_k)  + \lambda_k \e,  \\[1ex] 
         q_{k+1} = q_k + \alpha_k (\sbo_k-q_k) + \alpha_k \beta_k \Pi_M(p_k),  
    \end{array} 
\end{equation}
where 
\begin{equation}\label{eq:QPL-def0}
    \begin{array}{l}
        z = c - \frac{\e\tr c}{nm} \e  \in \R^{nm},\\[1ex]
        \wh{g}_k = q_k - \sbo_k \in \R^{nm},\quad 
        y_k = \wh{g}_k - \frac{\e\tr \wh{g}_k}{nm} \e \in \R^{nm},\\[1ex] 
        \delta_k = \left\{ \begin{array}{ll}
             0 & \text{if}\ q_k\tr (y_k+z) = 0, \\
            \frac{q_k\tr y_k}{q_k\tr (y_k+z)}  & \text{otherwise}
         \end{array} 
         \right.\in \R,\\[3ex]
        \lambda_k = \frac{\gamma}{nm(1-\gamma)} \e\tr \big((\delta_k-1) \wh{g}_k + \delta_k c\big) \in \R.
    \end{array}
\end{equation}  
\end{subequations}
Observe that in this case QPL uses the two additional vectors $c - \sbo_k$ and the all-one vector~$\e$ with adaptive coefficients in its update rule. 
The following theorem summarizes properties of the proposed QPL algorithm (see Section~\ref{appendix:proof thm QPL} for the proof).

\begin{Thm}[QPL] \label{thm:QPL}
Assume the learning rates $\alpha_k$ and $\beta_k$ are such that $\sum_k \alpha_k =  \infty$, $\sum_k \alpha_k^2  < \infty$, and $\beta_k \ra 0$. 
Then, the iterates~$q_k$ of QPL algorithm~\eqref{eq:QPL update gen} converge to $q\opt$ almost surely. 
Moreover, the algorithm has a per-iteration complexity of $\ord(nm^2)$.
\end{Thm} 

Regarding the preceding result, we note that the per-iteration time complexity of QPL is the same as that of the (synchronous) QL algorithm. 
We finish with the following remark on the \emph{asynchronous} implementation of QPL. 

\begin{Rem}[Asynchronous QPL]\label{rem:QPL asynch} 
The proposed QPL update rule~\eqref{eq:QPL update gen} can also be implemented in an asynchronous fashion. 
To be precise, this requires forming the approximate state-action transition matrix~$\wt{P}_k$ based on a single sample~$(s_k,a_k,\hat{s}_k^+, c(s_k,a_k))$ as follows
\begin{equation*}
\begin{array}{rl}
    \wt{P}_k = &\argmin\limits_{P \in \R^{nm \times nm}}  \Vert P - P\pr \Vert^2_F \\
               & \ \ \text{\emph{s.t.}} \ \   P \e = \e,\
                e_{(s_k,a_k)}\tr P q_k = {\gamma}^{-1} \big(  \sbo_k(s_k,a_k) -  c(s_k,a_k) \big).
\end{array}
\end{equation*}
Cf.~approximation~\eqref{eq:QPL-approx}.  
In particular, by using the uniform prior~$P\pr = \frac{1}{nm} E$, the corresponding update rule reads as
\begin{equation*}
\begin{array}{l}
    p_k =  \tau_k \big( \delta_k e_{(s_k,a_k)} +  \frac{\gamma(1+\delta_k)}{nm(1-\gamma)} \e \big),  \\[1ex] 
     q_{k+1} =  q_k + \alpha_k  \tau_k  e_{(s_k,a_k)} + \alpha_k \beta_k \Pi_M(p_k), 
\end{array}
\end{equation*}
where the scalar coefficients are given by 
\begin{equation*}
    \begin{array}{l}
        \tau_k = \sbo_k(s_k,a_k) - q_k(s_k,a_k), \quad 
        \rho_k = \frac{1}{nm} \e\tr q_k,\\[1ex] 
        \lambda_k = \big( \sbo_k(s_k,a_k) - c(s_k,a_k)  - \gamma \rho_k \big) \big( q_k(s_k,a_k) - \rho_k \big), \\[1ex] 
        \eta_k = \norm{q_k}_2^2 - \rho_k^2 -  \lambda_k, \quad 
        \delta_k = \left\{ \begin{array}{ll}
             0 & \text{if}\  \eta_k = 0, \\
            \lambda_k / \eta_k  & \text{otherwise}.
         \end{array} 
         \right.
    \end{array}
\end{equation*}  
Note that the preceding update rule leads to an update in all entries of the Q-function~$q_k$ in each iteration.
\end{Rem}

\section{Technical Proofs}\label{app:proofs}

\subsection{Proof of Theorem~\ref{thm:QPI}}\label{appendix:proof thm}

% \begin{proof}
First, let us show that the two equality constraints in the minimization problem~\eqref{eq:QPI approx} are linearly dependent if and only if $ u_k\tr v_k = 0$. 
In this regard, observe that the two equality constraints are linearly dependent if and only if $v_k = \rho \e$ for some $\rho \in \R$: 
For $\rho = 0$, the second equality constraint becomes trivial; and, for $\rho \neq 0$, the two constraints become equivalent. 
On the other hand, we have
\begin{align*}
    u_k\tr v_k = \big( v_k  - \frac{\e\tr v_k }{n} \e \big)\tr  v_k = v_k\tr (I - \frac{1}{n}E) v_k.
\end{align*} 
Then, since $I - \frac{1}{n}E$ is positive semi-definite with one zero eigenvalue corresponding to the eigenvector~$\e$, we have $u_k\tr v_k = 0$ if and only if $v_k = \rho \e$ for some $\rho \in \R$. 
Hence, the constraints in~\eqref{eq:QPI approx} are linearly dependent if and only if $u_k\tr v_k = 0$. 

We first consider the update rule~\eqref{eq:QPI update gen} for the case $u_k\tr v_k \neq 0$. 
Define $R \Let (\e, v_k),\; B \Let \big(\e, \gamma^{-1}(\bo_k -  c_k)\big) \in \R^{n\times 2}$ so that the minimization problem~\eqref{eq:QPI approx} can be written as
\begin{equation*}
    \wt{P}_k = \argmin_{P \in \R^{n\times n}} \left\{ \Vert P - P\pr\Vert^2_F\ :\ P R = B \right\}.
\end{equation*}
Note that, since $u_k\tr v_k \neq 0$ and hence $v_k$ and $\e$ are linearly independent, $R$ is of full column rank. The solution to the preceding problem is given by
\begin{equation*}
    \wt{P}_k = P\pr + (B - P\pr R) (R\tr R)^{-1} R\tr.
\end{equation*}
Now, observe that
\begin{align*}
    (B - P\pr R) 
    &= \left[\begin{array}{cc} \e  & \frac{1}{\gamma}(\bo_k -  c_k) \end{array}\right] - P\pr \left[\begin{array}{cc} \e  & v_k \end{array}\right] 
    = \left[\begin{array}{cc} \e-P\pr \e  & \frac{1}{\gamma}(\bo_k -  c_k) - P\pr v_k \end{array}\right] \\
    &= \left[\begin{array}{cc} \z  & \frac{1}{\gamma}(\bo_k -  c_k) - P\pr v_k \end{array}\right] 
    = \left[\begin{array}{cc} \z  & \gamma^{-1}w_k \end{array}\right], 
\end{align*}
where we used the assumption $P\pr \e = \e$. Also, 
\begin{align*}
    (R\tr R)^{-1} 
    &= \left( \left[\begin{array}{c} \e\tr \\ v_k\tr \end{array}\right] \left[\begin{array}{cc} \e  & v_k \end{array}\right] \right)^{-1} = \left[\begin{array}{cc} n & v_k\tr \e \\  v_k\tr \e & v_k\tr v_k \end{array}\right]^{-1} 
    = \frac{1}{n (u_k\tr v_k)} \left[\begin{array}{cc} v_k\tr v_k &  -\e\tr v_k \\  -\e\tr v_k & n \end{array}\right].
\end{align*}
and
\begin{align*}
    (R\tr R)^{-1}  R\tr
    &= \frac{1}{u_k\tr v_k} \left[\begin{array}{c} * \\  u_k\tr \end{array}\right].
\end{align*}
Therefore, we have
\begin{align*}
    \wt{P}_k =  P\pr + \gamma^{-1} (u_k\tr v_k)^{-1} w_k u_k\tr.
\end{align*}
For the case $u_k\tr v_k = 0$, as we discussed in the beginning of the proof, the two equality constraints in the minimization problem~\eqref{eq:QPI approx} become linearly dependent, and, in particular, the second constraint can be discarded. 
The solution to the problem~\eqref{eq:QPI approx} in this case is then $\wt{P}_k = P\pr$. 
Hence, the approximation~\eqref{eq:QPI approx} can be in general written as
\begin{equation}\label{eq:QPI approx sol}
    \wt{P}_k =  P\pr + \gamma^{-1} \tau_k w_k u_k\tr,
\end{equation}
where
\begin{equation*}
    \tau_k = \left\{ \begin{array}{ll}
    0 & \text{if}\ u_k\tr v_k = 0, \\
    (u_k\tr v_k)^{-1}  & \text{otherwise.}
    \end{array}\right.
\end{equation*}
That is, the approximation $\wt{P}_k$ is a rank-one update of the prior. 
Then, if $\tau_k = 0 $, we clearly have 
\[
\wt{G}_k = (I-\gamma \wt{P}_k)^{-1} = (I-\gamma P\pr)^{-1} = G\pr,
\] 
and, for the case $\tau_k \neq 0 $, we can use the the Woodbury formula to write 
\begin{align*}
    \wt{G}_k &= (G\pr^{-1} - \tau_k w_k u_k\tr )^{-1} 
    = G\pr + \frac{1}{\tau_k^{-1} -  u_k\tr (G\pr w_k)} (G\pr w_k) (G\pr\tr u_k)\tr  
    = G\pr + \eta_k \check{w}_k \check{u}_k\tr. 
\end{align*}
What remains to be shown is that $\eta_k$ is well-defined for $u_k\tr v_k \neq 0 $ (i.e., $\tau_k \neq 0$). 
First, we use $\bo_k = c_k + \gamma P_k v_k$ to write
\begin{align*}
    \eta_k^{-1} &= u_k\tr (v_k - \check{w}_k) 
    = u_k\tr \big(v_k - G\pr (\bo_k - c_k - \gamma P\pr v_k)\big) \\
    &= u_k\tr \big(v_k - G\pr (\gamma P_k v_k - \gamma P\pr v_k)\big),
\end{align*}
Next, since $P\pr = \gamma^{-1}(I - G\pr^{-1})$ and using the fact that 
$v_k = u_k + \frac{\e\tr v_k}{n} \e$, we have
\begin{align*}
    \eta_k^{-1} &= u_k\tr \big(I - \gamma G\pr ( P_k - P\pr) \big)v_k 
    = u_k\tr G\pr (I-\gamma P_k) v_k \\ 
    &= u_k\tr G\pr (I-\gamma P_k) (u_k + \alpha_k \e) 
    = u_k\tr G\pr (I-\gamma P_k) u_k + \frac{\e\tr v_k}{n} u_k\tr G\pr (I-\gamma P_k) \e.
\end{align*}
Now, note that $P_k \e = \e$ and $G\pr \e =  (1-\gamma)^{-1} \e$. Hence,
\begin{align*}
    \eta_k^{-1} &= u_k\tr G\pr (I-\gamma P_k) u_k + \frac{\e\tr v_k}{n} u_k\tr  \e
    = u_k\tr G\pr (I-\gamma P_k) u_k,
\end{align*}
where we also used the fact that $u_k\tr  \e = 0$. 
Then, since the matrices $G\pr$ and $(I-\gamma P_k)$ are non-singular with their eigenvalues having strictly positive real parts (because the eigenvalues of $P_k$ all reside within the unit disc), we have
\begin{align*}
    u_k\tr G\pr (I-\gamma P_k) u_k = 0 &\iff u_k = \z 
    \iff  v_k = \rho \e \ \ \text{for some} \ \ \rho \in \R 
    \iff  u_k\tr v_k =0.
\end{align*}
That is, $\eta_k$ is well-defined for $u_k\tr v_k \neq 0 $.

Next, we consider the rate of convergence of the safeguarded QPI update rule.
We will use induction to show that 
\begin{equation}\label{eq:proof_induction}
\theta_k \leq \gamma^{k} \theta_0, \quad \forall k\geq 0.
\end{equation}
The base case $k=0$ holds trivially. 
Let us now assume the inequality \eqref{eq:proof_induction} holds for some $k\ge0$ and recall that $\bo$ is a $\gamma$-contraction in $\infty$-norm. 
Then, the safeguarding against standard VI as in \eqref{eq:QPI safegaurd} implies that 
\begin{align*}
        \theta_{k+1}  &\leq \max \{ \gamma^{k+1} \theta_{0},\; \gamma \theta_{k} \} \leq  \gamma^{k+1} \theta_{0},
\end{align*}
where we used the induction hypothesis for the last inequality. This completes the proof.
    
Finally, the per-iteration time complexity of each iteration of the safeguarded QPI update rule: 
The update rule~\eqref{eq:QPI update} requires $\ord(n^2m)$ operations for computing the vectors $\bo_k = \bo (v_{k})$, $\ord(n^2)$ operations for the matrix-vector multiplication, and $\ord(n)$ operations for the vector additions. 
Computing the objects in \eqref{eq:QPI-def0} involves vector/matrix additions and matrix-vector multiplications (all of size $n$) and hence requires $\ord(n^2)$ operations.  
For the safeguarding~\eqref{eq:QPI safegaurd}, we need to compute $\bo_{k+1} = \bo (v_{k+1})$ which again requires $\ord(n^2m)$ operations. 
Summing up the aforementioned complexities, we derive the total time complexity to be $\ord(n^2m)$.
% \end{proof}

\subsection{Proof of Corollary~\ref{cor:QPI-A}}\label{appendix:proof cor}

% \begin{proof}
The result follows from Theorem~\ref{thm:QPI} by plugging in $P\pr = \frac{1}{n}E$ and $G\pr = I+\frac{\gamma}{n(1-\gamma)}E$ and simplifying the expression. 
In particular, we note that the condition $v_k\tr (y_k+z_k) = 0$ in~\eqref{eq:QPI-A-def0} is equivalent to the condition $u_k\tr v_k = 0$ in \eqref{eq:QPI-def0}. 
To see this, recall that $u_k\tr v_k = 0$ if and only if $v_k = \rho \e$ for some $\rho \in \R$; see the first part of the proof of Theorem~\ref{thm:QPI} in Section~\ref{appendix:proof thm}. 
Also, observe that
\begin{align*}
    v_k\tr (y_k+z_k) &= v_k\tr \big( g_k+c_k  - \frac{\e\tr (g_k+c_k) }{n} \e \big) 
    = v_k\tr \big( v_k - \bo_k+c_k  - \frac{\e\tr (v_k - \bo_k+c_k) }{n} \e \big) \\
    &= v_k\tr (I - \frac{1}{n}E) (v_k - \bo_k+c_k) 
    = v_k\tr (I - \frac{1}{n}E)(I - \gamma P_k) v_k,
\end{align*}
where, for the last equality, we used $\bo_k = c_k + \gamma P_k v_k$. 
Then, since $u_k = (I - \frac{1}{n}E)v_k = v_k - \frac{\e\tr v_k}{n} \e $, we have
\begin{align*}
    v_k\tr (y_k+z_k) &= u_k\tr (I - \gamma P_k) (u_k + \frac{\e\tr v_k}{n} \e) 
    = u_k\tr (I - \gamma P_k) u_k + \frac{\e\tr v_k}{n}  u_k\tr (I - \gamma P_k) \e \\
    &= u_k\tr (I - \gamma P_k) u_k + \frac{\e\tr v_k}{n} (1-\gamma)  u_k\tr \e 
    = u_k\tr (I - \gamma P_k) u_k,
\end{align*}
where, for the last equality, we used the fact that $u_k\tr \e = 0$. 
Finally, since $(I - \gamma P_k)$ is non-singular with its eigenvalues having strictly positive real parts (because the eigenvalues of $P_k$ all reside within the unit disc), we have 
\begin{align*}
    v_k\tr (y_k+z_k) = 0  &\iff u_k = \z 
    \iff  v_k = \rho \e \ \ \text{for some} \ \ \rho \in \R.
\end{align*}
This completes the proof. 
% \end{proof}

\subsection{Proof of Lemma~\ref{lem:backtrack}}\label{appendix:proof lem backtrack} 

Recall that $T$ is a $\gamma$-contraction observe that 
\begin{align*}
    \theta_{\bo(v)} &= \norm{\bo(v)-\bo(\bo(v))}_{\infty}  
    \leq \gamma \norm{v-\bo(v)}_{\infty} = \gamma \theta_v.
\end{align*}
Also, note that the map $\theta_{(\cdot)}$ is $(1+\gamma)$-Lipschitz continuous: For any $w_1,w_2 \in \R^n$, we have
\begin{align*}
     |\theta_{w_1}-\theta_{w_2}| &= |\norm{w_1-\bo(w_1)}_{\infty} - \norm{w_2 - \bo(w_2)}_{\infty}|  
     \leq \norm{w_1-\bo(w_1) - w_2 + \bo(w_2)}_{\infty}\\
     &\leq \norm{w_1-w_2}_{\infty} + \norm{\bo(w_1)-\bo(w_2)}_{\infty}
     \leq (1+\gamma) \norm{w_1-w_2}_{\infty}.
\end{align*}
Setting $w_1 = w$ with $\norm{w-\bo(v)}_{\infty} \leq \frac{\gamma'-\gamma}{1+\gamma}\theta_v$ and $w_2 = \bo(v)$, we then have
\begin{align*}
    |\theta_{w}-\theta_{\bo(v)}| &\leq (1+\gamma) \norm{w-\bo(v)}_{\infty} \leq (\gamma'-\gamma)\theta_v.
\end{align*}
Hence,  
\begin{align*}
    \theta_{w} & \leq \theta_{\bo(v)} + (\gamma'-\gamma)\theta_v \leq \gamma'\theta_v.
\end{align*}
This completes the proof. 

\subsection{Proof of Theorem~\ref{thm:QPI-B superlinear}}\label{appendix:proof thm QPI-B superlinear}

In what follows, we show that the proposed algorithm is an instance of semismooth QNM with a structured least-change secant update and a standard backtracking scheme. 
To this end, observe that we are looking for the unique root~$v\opt$ of the map $g(v) = v - \bo(v)$. 
First, note that $g$ is piece-wise affine, Lipschitz-continuous, and, in particular, \emph{strongly semismooth}~\cite[Prop.~7.4.7]{facchinei2003finite}. 
Moreover, the any generalized Jacobian $H_k \in \partial g(v) \subset \{I - \gamma P \;:\; P\ge 0,\ P\e = \e\}$ is \emph{non-singular} at any point~$v$ (note that for a row stochastic matrix $P$ and $\gamma \in (0,1)$, we have $(I - \gamma P)^{-1} = \sum_{t=0}^{\infty} \gamma^tP^t$). 

Denote $g_k \Let g(v_k)$ and $\wt{H}_k \Let \wt{G}_k^{-1}$ for each $k\geq0$. 
The update rule~\eqref{eq:QPI backtrack} then can be equivalently written as
\begin{equation*}\label{eq:QPI backtrack copy}
     v_{k+1} = v_k - \alpha_k \wt{H}_k^{-1} g_k - (1-\alpha_k) g_k,     
\end{equation*}
where $\alpha_k$ is determined using the backtracking scheme~\eqref{eq:QPI backtrack 1} ensuring
\[
\norm{g_{k+1}}_{\infty} \leq \gamma' \norm{g_{k}}_{\infty}, 
\]
for some $\gamma' \in (\gamma, 1)$. 
This means that, as long as the approximations~$\wt{H}_k$ of the Jacobian of $g$ are uniformly bounded, we have the global linear convergence of the iterates $v_k$ to $v\opt$ with rate~$\gamma'$; see Lemma~\ref{lem:backtrack}. 

Finally, let us look at the approximations~$\wt{H}_k$ of the Jacobian of $g$. Using the change of variable $\wt{P}_k = \gamma^{-1}(I-\wt{H}_k)$, we can show that
\begin{equation*}\label{eq:QPI-B gain selection}
\begin{array}{l}
    \text{{\bf if} $\pi_{v_k} = \pi_{v_{k-1}}$}  \\
    \hspace{1cm}\begin{array}{rl}
        \wt{H}_k = \argmin\limits_{H \in \R^{n \times n}} & \Vert H - \wt{H}_{k-1} \Vert^2_F \\
        \text{s.t.} & H \e = (1-\gamma)\e,\  
          H v_k = g_k + c_k , \ 
         H (v_k-v_{k-1}) = g_k-g_{k-1};
        \end{array} \\
    \text{{\bf else}}  \\
    \hspace{1cm}\begin{array}{rl}
        \wt{H}_k = \argmin\limits_{H \in \R^{n \times n}} & \Vert H - \wt{H}_{k-1} \Vert^2_F \\
        \text{s.t.} & H \e = (1-\gamma)\e, \
         H (v_k-v_{k-1}) = g_k-g_{k-1}.
    \end{array}
\end{array}
\end{equation*}
One can now clearly observe that the approximation $\wt{H}_k$ satisfies the \emph{secant} condition, i.e.,  
\[
\wt{H}_k (v_k-v_{k-1}) = g_k-g_{k-1},
\]
besides the structural linear constraints that the true (generalized) Jacobian must satisfy. 
This implies that the proposed algorithm indeed involves a least-change in $\Vert \wt{H}_{k} - \wt{H}_{k-1} \Vert_F$ subject to secant and structural constraints. 
Now, since $g$ is strongly semismooth, there exists $H_k \in \partial g(v_k)$ such that~\cite[Def.~7.4.2]{facchinei2003finite}
\[
\|g_{k} - g_{k-1} - H_k (v_{k}-v_{k-1})\| = o(\|v_{k}-v_{k-1}\|).
\]
Therefore,
\[
\| (\wt{H}_k - H_k) (v_{k}-v_{k-1}) \| = o(\|v_{k}-v_{k-1}\|).
\]
That is, Dennis-Mor\'{e} condition is satisfied and, by \cite[Thm.~4.2]{sun1997newton}, the convergence is locally superlinear.   

\subsection{Proof of Theorem~\ref{thm:QPL}}\label{appendix:proof thm QPL}

% \begin{proof}
The per-iteration time complexity of each iteration of the QPL update rule is as follows: 
The update rule~\eqref{eq:QPL update gen} requires $\ord(nm^2)$ operations for computing the vectors $\sbo_k = \sbo (q_k)$, and $\ord (nm)$ operations for computing the step-sizes $\delta_k$ and $\lambda_k$ and the vector additions.  
Summing up the aforementioned complexities, the total time complexity is $\ord(nm^2)$. 

Regarding the convergence, define 
\begin{align*}
    F(q) \Let \bar{\bo}(q) - q,\ 
    V_k \Let \sbo_k - \bo(q),\
    d_k \Let \beta_k \Pi_M(p_k),
\end{align*}
where $\bar{\bo}$ is the corresponding true Bellman operator for Q-functions, i.e., $\bar{\bo}(q) = \EE[\sbo (q)]$ for each $q$. 
Then, the convergence of $q_k$ to the unique equilibrium $q\opt$ of the map $F$ follows from \cite[Cor.~3.2]{perkins2013asynchronous}. 
In particular, observe that $d_k \ra 0$ since $\beta_k \ra 0$ and $\norm{\Pi_M(p_k)}_\infty \leq M$. 
 
% \end{proof}

\section{Numerical Simulations}\label{sec:experiments} 
We now compare the performance of the proposed algorithms with that of the standard existing algorithms for the optimal control of different MDPs. 
See Appendix~\ref{appendix:MDPs} for a description of the considered MDPs. 
To this end, we first focus on the proposed algorithms with uniform priors corresponding to update rules~\eqref{eq:QPI-A update gen} and \eqref{eq:QPL update gen} in Sections~\ref{sec:exp-mb} and \ref{sec:exp-mf}, respectively. 
The results of numerical experiments with alternative priors are then reported in Section~\ref{app:prior}.

\subsection{Model-based algorithms}\label{sec:exp-mb}

\begin{figure*}
\centering
\subfloat[]{{\includegraphics[clip, trim=2cm 0cm 2cm 0cm,width=1\linewidth]{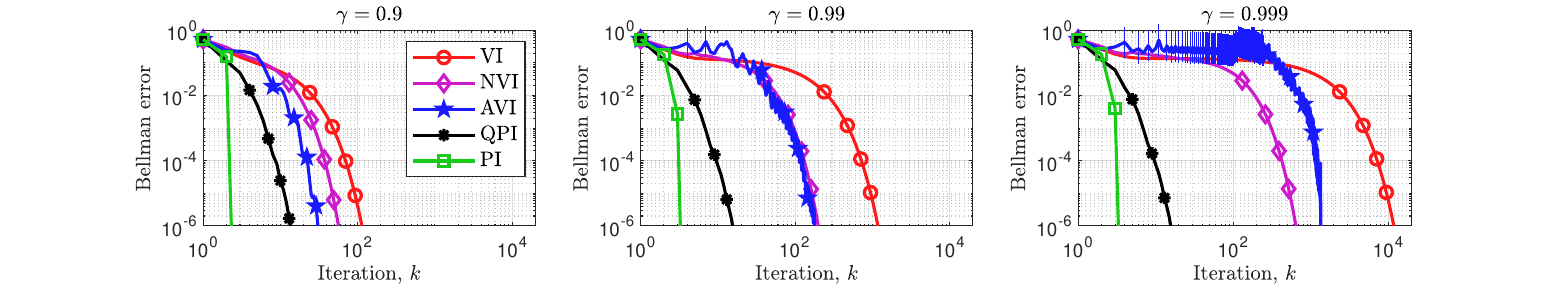}}\label{fig:mb_garnet}}

\subfloat[]{{\includegraphics[clip, trim=2cm 0cm 2cm 0cm,width=1\linewidth]{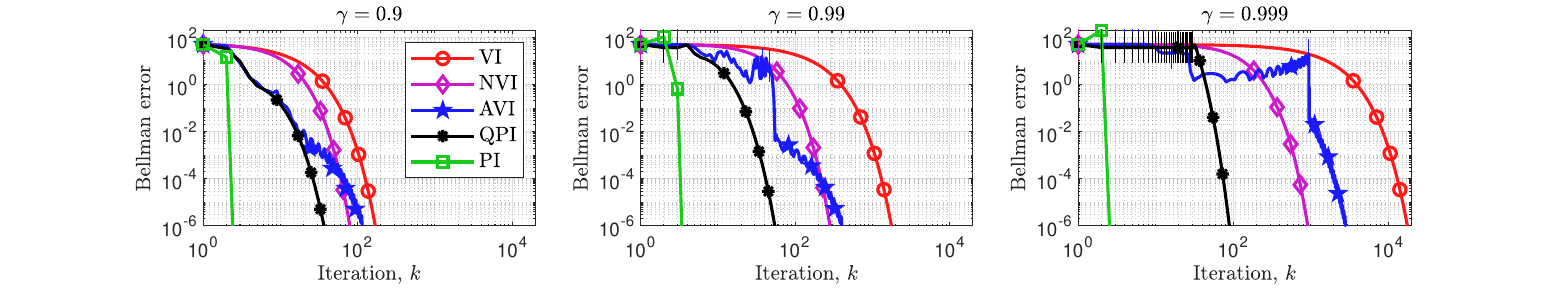}}\label{fig:mb_health}}

\caption{Performance of model-based algorithms for three values of $\gamma$: (a)~Garnet MDP; (b)~Healthcare MDP. The bars indicate the iterations at which the safeguard is activated (for NVI, AVI, and QPI).}
\label{fig:mb_instance} 
\end{figure*}

\begin{figure}
\centering
\includegraphics[clip, trim=0cm 0cm 0cm 0cm,width=.29\linewidth]{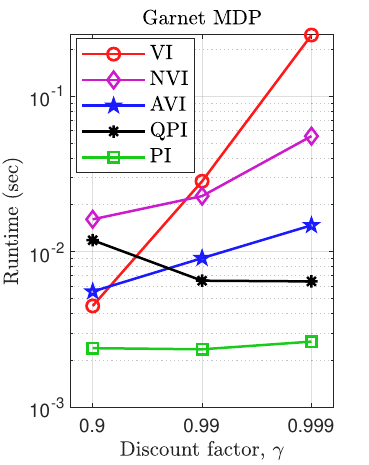}
\includegraphics[clip, trim=0cm 0cm 0cm 0cm,width=.29\linewidth]{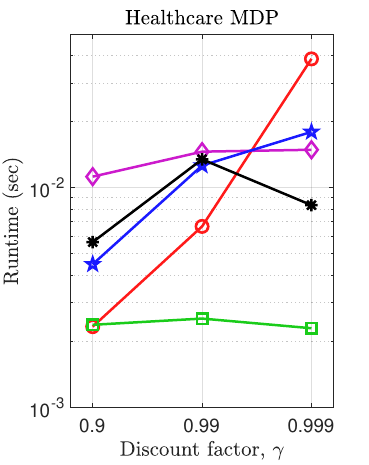}
\caption{The running time of the model-based algorithms for three values of $\gamma$ corresponding to Figure~\ref{fig:mb_instance}.}
\label{fig:mb_rt}
\end{figure}

For model-based algorithms we consider two MDPs: a randomly generated Garnet MDP~\cite{archibald1995generation} and the Healthcare MDP~\cite{goyal2019first} with an \emph{absorbing} state. 
The proposed QPI algorithm~\eqref{eq:QPI-A update gen} is compared with the following algorithms: VI (value iteration); NVI (VI with Nesterov acceleration)~\cite{goyal2019first}; AVI (VI with Anderson acceleration)~\cite{geist2018anderson}; and, PI (policy iteration). 
For AVI, we use a memory of one leading to a rank-one update (of the identity matrix) for approximating the Hessian so that it is comparable with the rank-one update of the uniform distribution for approximating the transition matrix in QPI. 
See Appendix~\ref{app:NVI and AVI} for the exact update rules of NVI and AVI. 
We note that since NVI and AVI are not guaranteed to converge, we safeguard them using VI (using the same safeguarding rule~\eqref{eq:QPI safegaurd} used for QPI). 
All the algorithms are initialized by $v_0 = \z$ with termination condition $\norm{v_k - \bo (v_k)}_{\infty} \leq 10^{-6}$. 
The results of the simulations are provided in Figures~\ref{fig:mb_instance} and \ref{fig:mb_rt}.

In Figure~\ref{fig:mb_instance}, VI, NVI, and AVI show a linear convergence with a rate depending on~$\gamma$ in both MDPs. 
In particular, as we increase $\gamma$ from $0.9$ to $0.999$, we observe more than a tenfold increase in the number of iterations required for these algorithms to terminate. 
This is expected since these algorithms only use first-order information and their convergence rate is determined by $\gamma$.  

Figure~\ref{fig:mb_instance} also shows that for both MDPs, PI converges with a quadratic rate in 3 to 5 iterations, independent of $\gamma$. 
Now, observe that QPI is the only algorithm showing a similar behavior as PI and terminating in approximately the same number of iterations, independent of $\gamma$, in both MDPs. 
Moreover, comparing the performance of QPI with AVI (its counterpart in the class of QNMs), we also see the importance of newly introduced linear constraints and prior in the approximation of the transition matrix. 
Moreover, observe that QPI's safeguard is activated for Healthcare MDP as shown in Figure~\ref{fig:mb_health} (one instance for $\gamma=0.99$ and multiple instances for $\gamma=0.999$). 
In this regard, we note that for QPI, we have observed that the activation of the safeguard is particularly due to the existence of \emph{absorbing states} in the MDP as is the case for the Healthcare MDP. 

Figure~\ref{fig:mb_rt} reports the corresponding running times of the algorithms. 
The reported running times are in line with convergence behaviors seen in Figure~\ref{fig:mb_instance} and the theoretical time complexity of these algorithms. 
In particular, PI and QPI are the only algorithms with running time less sensitive to $\gamma$ for both of the considered MDPs. 
In this regard, we note that since the size of the MDPs considered in our numerical simulations is relatively small, PI is the fastest algorithm despite the fact that it requires a matrix inversion. 

\begin{figure*}
\centering
\subfloat[]{{\includegraphics[clip, trim=1.5cm 0cm 1.5cm 0cm,width=1\linewidth]{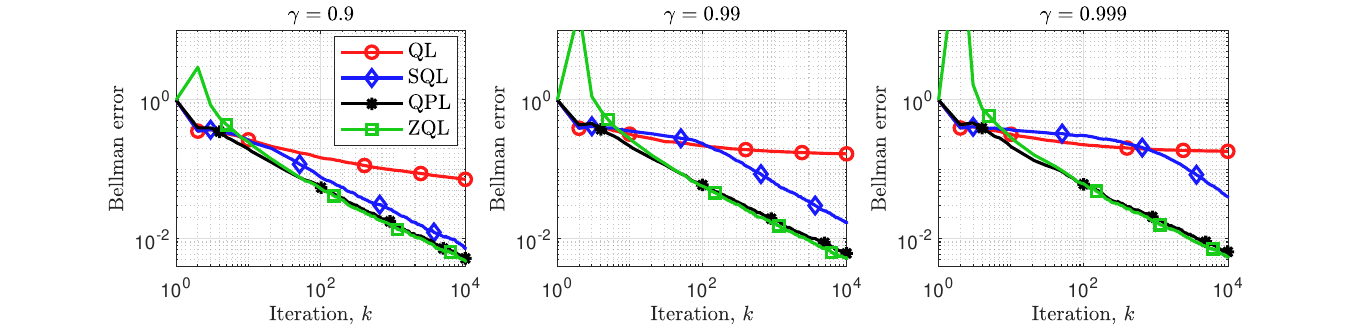}}\label{fig:mf_garnet}}

\subfloat[]{{\includegraphics[clip, trim=1.5cm 0cm 1.5cm 0cm,width=1\linewidth]{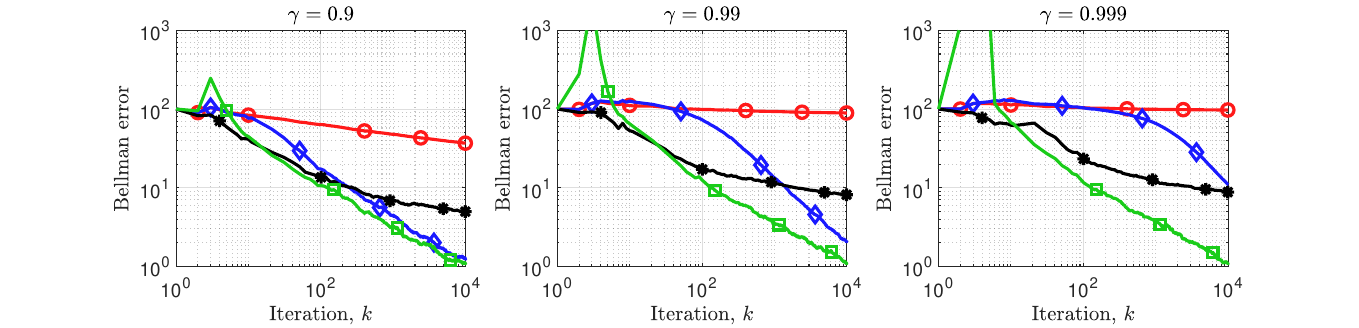}}\label{fig:mf_graph}}

\caption{Performance of model-free algorithms (averaged over 20 runs) for three values of $\gamma$: (a)~Garnet MDP; (b)~Graph MDP.}
\label{fig:mf_instance} 
\end{figure*}

\subsection{Model-free algorithms}\label{sec:exp-mf}

For model-free algorithms we also consider two MDPs: again a randomly generated Garnet MDP~\cite{archibald1995generation} and 
the Graph MDP~\cite{devraj2017zap}. 
The proposed QPL algorithm~\eqref{eq:QPL update gen} with $\alpha_k = \frac{1}{(k+1)}$ and $\beta_k = \frac{1}{(k+1)^{0.1}}$ is compared with the following algorithms: QL (Q-learning) as in~\eqref{eq:QL} with $\alpha_k = \frac{1}{(k+1)}$; 
SQL (speedy QL)~\cite{ghavamzadeh2011speedy}; 
and, ZQL (zap QL)~\cite{devraj2017zap}. 
See Appendix~\ref{app:NVI and AVI} for the exact update rules of SQL and ZQL. 
All the algorithms are initialized by $q_0 = \z_{nm}$ and terminated after $K = 10^4$ iterations with a synchronous sampling of all state-action pairs at each iteration. 
% For QPL, the safeguard is activated after $K_{\text{sg}} = 100$ iterations. 
For each algorithm, we report the \emph{average} of the Bellman error $\norm{q_k - \bo (q_k)}_{\infty}$ over 20 runs of the algorithm. 
The results of the simulations are provided in Figure~\ref{fig:mf_instance} and Table~\ref{tab:mf_rt}. 

\begin{table}
\caption{The running time (in seconds) of the model-free algorithms over $K=10^4$ iterations (averaged over 20 runs) for $\gamma = 0.9$ corresponding to Figure~\ref{fig:mf_instance}.}
\label{tab:mf_rt}
%\vskip 0.15in
\begin{center}
\begin{small}
%\begin{sc}
\begin{tabular}{l|cccc|c}
\toprule
 & QL & SQL & QPL & ZQL & Sampling \\
\midrule
Garnet & $1.7$ & $3.1$ & $3.3$ & $16$ & $71$ \\
Graph  & $0.16$ & $0.27$ & $0.31$ & $0.59$ & $6.5$ \\
\bottomrule
\end{tabular}
%\end{sc}
\end{small}
\end{center}
%\vskip -0.1in
\end{table}

As can be seen in Figures~\ref{fig:mf_garnet} and \ref{fig:mf_graph}, the performance of QL and SQL (the first-order methods) deteriorates as $\gamma$ increases for both MDPs. 
However, for these MDPs, ZQL (the second-order method that estimates the transition matrix by averaging over the samples) leads to almost the same error level after a fixed number of iterations for different values of $\gamma$. 

Figure~\ref{fig:mf_instance} shows that the performance of QPL is not as consistent as its model-based counterpart: 
QPL has the same rate of convergence as ZQL for Garnet MDP (Figure~\ref{fig:mf_garnet}), while it deviates from ZQL as $k$ increases for Graph MDP (Figures~\ref{fig:mf_graph}). 
This means that for structured MDPs, QPL may not lead to a better performance compared to SQL or ZQL. 
In this regard, we note that the model-free QPL algorithm uses an approximation of the transition matrix, which is constructed based on sampled data; see \eqref{eq:QPL-approx}. 
This use of sampling on top of approximation can be the reason behind the poor performance of the model-free QPL algorithm for structured MDPs. 

Finally, we note that the running times reported in Table~\ref{tab:mf_rt} also align with the corresponding theoretical time complexities of these algorithms. 
In particular, QPL and SQL require almost the same amount of time, which is slightly more than QL and less than ZQL. (We report the runtime only for $\gamma = 0.9$ because it is independent of $\gamma$).  
Note that the time required for generating the samples is reported separately in Table~\ref{tab:mf_rt}, which is indeed the dominating factor in the actual runtime of the model-free algorithms.

\subsection{QPI and QPL with different priors}\label{app:prior}

\begin{figure*}
\centering
\subfloat[]{{\includegraphics[clip, trim=1.5cm 0cm 1.5cm 0cm,width=1\linewidth]{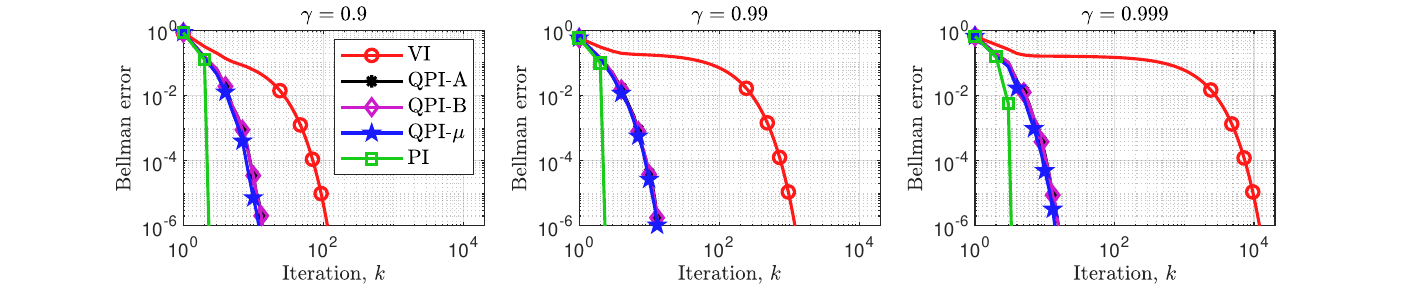}}\label{fig:mb_garnet_prior}}

\subfloat[]{{\includegraphics[clip, trim=1.5cm 0cm 1.5cm 0cm,width=1\linewidth]{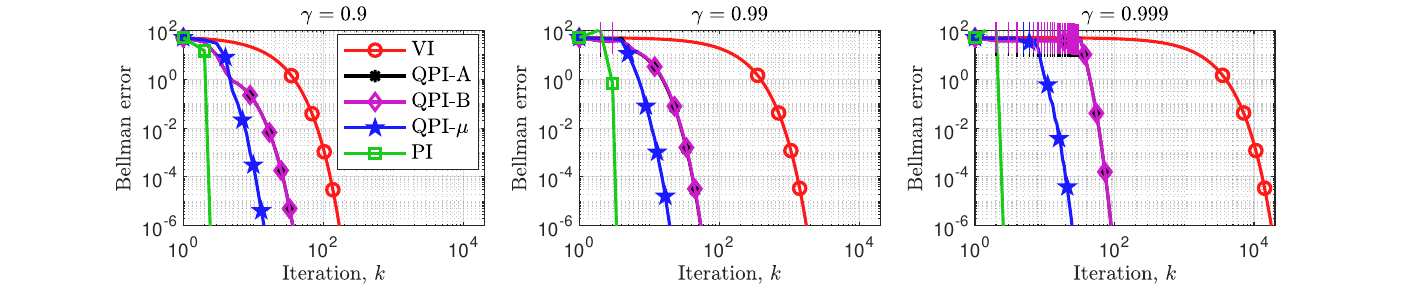}}\label{fig:mb_health_prior}}

\caption{Performance of model-based algorithms for three values of $\gamma$ and three different priors: (a)~Garnet MDP; (b)~Healthcare MDP. The bars indicate the iterations at which the safeguard is activated in QPI.}
\label{fig:mb_prior} 
\end{figure*}

\begin{figure*}
\centering
\subfloat[]{{\includegraphics[clip, trim=1.5cm 0cm 1.5cm 0cm,width=1\linewidth]{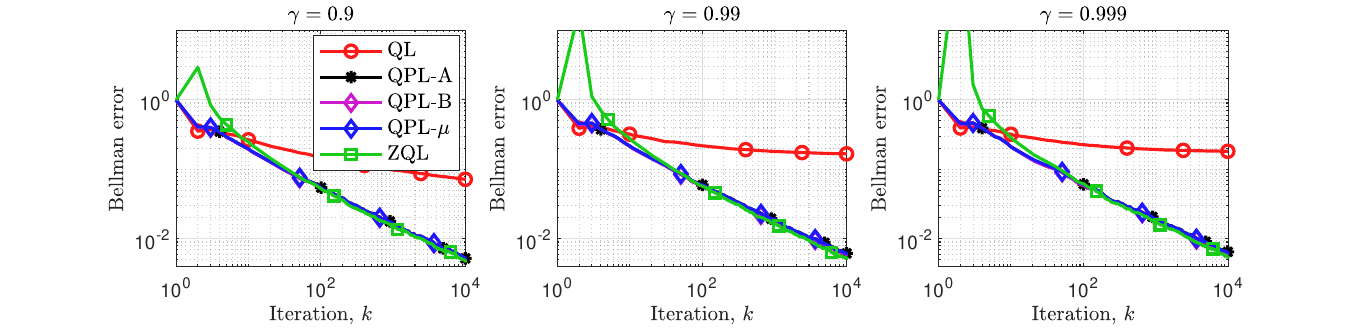}}\label{fig:mf_garnet_prior}}

\subfloat[]{{\includegraphics[clip, trim=1.5cm 0cm 1.5cm 0cm,width=1\linewidth]{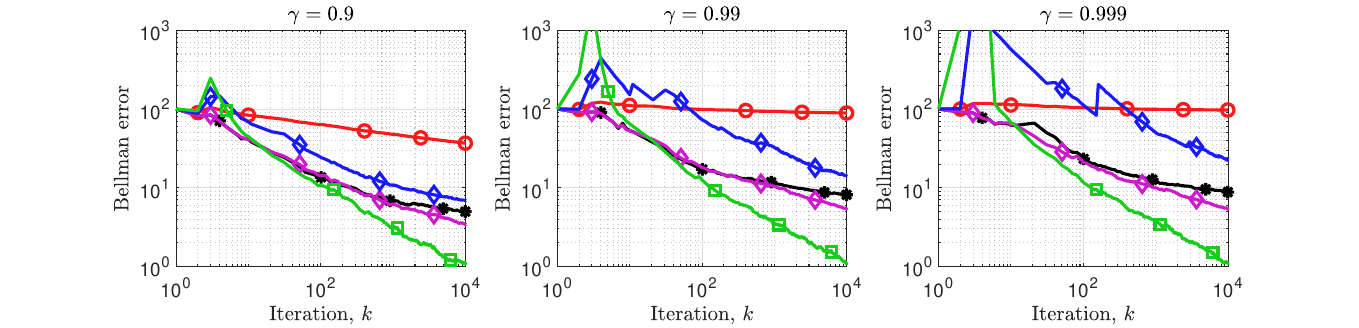}}\label{fig:mf_graph_prior}}

\caption{Performance of model-free algorithms (averaged over 20 runs) for three values of $\gamma$ and three different priors: (a)~Garnet MDP; (b)~Graph MDP.}
\label{fig:mf_prior} 
\end{figure*}

Figures~\ref{fig:mb_prior} and~\ref{fig:mf_prior} report the result of our numerical simulations for the QPI and QPL algorithms, respectively, with the three choices of the prior: (i)~QPI/L-A with a uniform prior $P\pr \propto \e\e\tr$, (ii)~QPI/L-B with recursive prior $P\pr = \wt{P}_{k-1}$, and (iii)~QPI/L-$\mu$ with prior $P\pr = P^{\mu}$ and $\mu$ being the stochastic policy choosing actions uniformly at random so that the prior has the same sparsity pattern as the true transition probability matrix. 

As depicted in Figures~\ref{fig:mb_garnet_prior} and \ref{fig:mf_garnet_prior}, the experiments with alternative priors shows no improvement in the performance of the QPI and QPL algorithms in comparison with the uniform prior for random Garnet MDPs. 
For structured MDPs, however, we observe contradictory results as shown in Figures~\ref{fig:mb_health_prior} and \ref{fig:mf_graph_prior}: 
Using a structured prior leads to a significant improvement in the performance of the (model-based) QPI algorithm for Healthcare MDP, while using a structured prior significantly deteriorates the performance of the (model-free) QPL algorithm for Graph MDP.  

\section{Limitations and Future Research}\label{sec:conclusion} 

In this paper, we proposed the model-based quasi-policy iteration (QPI) algorithm and its model-free counterpart, the quasi-policy learning (QPL) algorithm. 
The proposed algorithms were particularly inspired by the quasi-Newton methods and employed a novel approximation of the ``Hessian'' by using two new linear constraints specific to MDPs.

The main drawback of the proposed algorithms, similar to other accelerated VI schemes in the literature, is the need for safeguarding to ensure convergence. 
Our experiments in Section~\ref{sec:experiments} showed examples of MDPs in which the safeguard is activated. 
Our modified implementation of Section~\ref{sec:QPI backtrack} guarantees the linear convergence while removing the need for safeguarding and using backtracking instead. 
Another possible solution is the use of the operator splitting method introduced in~\cite{NEURIPS2022_fa809df3} for policy evaluation. 
In this regard, we note that the proposed QPI algorithm is essentially the PI algorithm in which the policy evaluation step uses the approximation~$\wt{P}_k$ in~\eqref{eq:QPI approx} instead of the true transition matrix~$P_k$ and the cost~$\wt{c}_k = c_k + \gamma(P_k - \wt{P}_k) v_k$ instead of the true cost~$c_k$. 
However, the convergence requires $\wt{P}_k$ to be close to $P_k$. 
To be precise, a sufficient condition is $\Vert P_k-\wt{P}_k\Vert _{\infty} \leq 1-\gamma$~\cite[Thm.~1]{NEURIPS2022_fa809df3}, which is difficult to achieve for a low-rank approximation~$\wt{P}_k$ of $P_k$.

Another limitation of the current work is the lack of a theoretical guarantee for the empirically observed improvement in the convergence rate, particularly for the model-based QPI algorithm. 
The main difficulty to be addressed is the fact that the linear constraints in~\eqref{eq:QPI approx} are not the standard secant conditions used in QNMs. 
The modified implementation of Section~\ref{sec:superlinear} partially addresses this issue by including the secant conditions in its approximation. 
An interesting result, however, is the establishment of a local superlinear convergence rate for QPI without imposing the secant condition. 
Another possibility is to use the results for Anderson acceleration in~\cite{Evans2020} to establish an improved linear rate for convergence. 
To that end, similar to what is done in~\cite{sun2021damped}, one needs to use a smoothed version of the Bellman operator, e.g., by replacing the \emph{max} operation with a \emph{soft-max} operation in the Bellman operator.  

The proposed algorithms in this study heavily rely on the approximation~\eqref{eq:QPI approx} of the transition matrix. 
As we discussed, this approximation easily allows for incorporation of different priors, e.g., a prior with the same sparsity pattern as the true transition matrix, or, the recursive prior. 
However, our numerical simulations using these alternative priors did not demonstrate a clear improvement in the performance of the proposed algorithms, highlighting the need for further investigation on other MDPs. 
In this regard, we also note that the main drawback of the approximation~\eqref{eq:QPI approx} is that it does not allow for a computationally efficient incorporation of other constraints, such as non-negativity constraints. 
A promising future research direction is the development of alternative approximation schemes that allow such constraints to be included at a reasonable computational cost. 

\appendix

\section{MDPs of the numerical simulations}\label{appendix:MDPs}
\textbf{Garnet MDP.} The considered Garnet MDP~\cite{archibald1995generation} is generated randomly with $n=50$ states, $m = 5$ actions, and the branching parameter~$n_b = 10$. 
For each state-action pair~$(s,a)$, we first form the set of reachable next states~$\{s_1^+,\ldots,s_{n_b}^+ \}$ chosen uniformly at random from the state space~$\{1,\ldots,n\}$. 
Then, the corresponding probabilities are formed by choosing the points $p_i\in [0,1], i=1,\ldots,n_b-1$, uniformly at random, and setting $\PP (s_i^+|s,a) = p_{i} - p_{i-1}$ with $p_0 = 0$ and $p_{n_b} = 1$. 
The stage cost $c(s,a)$ for each state-action pair~$(s,a)$ is also chosen uniformly at random from the interval $[0,1]$. 

\textbf{Healthcare MDP.} The considered Healthcare MDP is borrowed from~\cite{goyal2019first}. The MDP has 6 states corresponding to the deteriorating health condition of a patient with the last state $n=6$ being an absorbing state representing the mortality terminal state. For each of the first five states, one can choose three inputs $m\in\{1,2,3\}$ corresponding to increasing levels of drug dosage for treatment. The goal is to minimize the invasiveness of the treatment while avoiding the terminal state. 
For the transition probabilities, we refer the reader to~\cite[Fig.~D.1]{goyal2019first}. 
The cost function is chosen to be $c(n,m) = \sum_{n^+=1}^6 n^+ \PP(n^+|n,m) + m$ for each $n\in\{1,2,3,4,5\}$ and $m\in \{1,2,3\}$ and $c(6,1) = 50$. 

\textbf{Graph MDP.} The considered Graph MDP is borrowed from~\cite{devraj2017zap}. 
The MDP has 18 state-action pairs in total and corresponds to a simple path-finding problem. 
We refer the reader to~\cite[Sec.~3]{devraj2017zap} for the description of the MDP. 

\section{Accelerated VI and QL algorithms}\label{app:NVI and AVI}

The update rules are as follows: 

\noindent\textbf{(1)}~Nesterov accelerated VI~\cite{goyal2019first} -- with $v_{-1} = v_0$:
    \begin{align*}
        &y_k = v_k + \gamma^{-1} (1- \sqrt{1-\gamma^2}) ( v_k - v_{k-1}),\\
        &v_{k+1} = y_k - (1+\gamma)^{-1}  \big( y_k - \bo(y_k) \big).
    \end{align*}

\noindent\textbf{(2)}~Anderson accelerated VI~\cite{geist2018anderson} -- with $v_{-1} = v_0$:
    \begin{align*}
        &y_k = v_k - v_{k-1},\\
        &z_k = \bo(v_k) - \bo(v_{k-1}) ,\\
        &\delta_k = \frac{y_k\tr \big( v_k - \bo(v_k) \big)}{y_k\tr (y_k - z_k)}\ \text{if}\  y_k\tr (y_k - z_k) \neq 0,\ =0\ \text{o.w.}, \\
        &v_{k+1} = (1-\delta_k) \bo(v_k) + \delta_k \bo(v_{k-1}).
    \end{align*}

\noindent\textbf{(3)}~Speedy QL~\cite{ghavamzadeh2011speedy} -- with $q_{-1} = q_0$, $\alpha_k = \frac{1}{k+1}$:
    \begin{equation*}\label{eq:SQL}
    \left\{\begin{array}{l}
         d_{k} = \sbo (q_{k}) - \sbo (q_{k-1}),\\
         q_{k+1} = q_k -  \alpha_k \big(q_k -\sbo (q_{k-1})\big) + (1-\alpha_k)d_{k}. 
    \end{array}\right.
    \end{equation*}

\noindent\textbf{(4)}~Zap QL~\cite{devraj2017zap} -- with $D_{-1} = \z$, $\alpha_k = \beta_k = \frac{1}{k+1}$:
    \begin{equation*}\label{eq:ZQL}
    \left\{\begin{array}{l}
        D_k = (1-\beta_k) D_{k-1} + \beta_k  \big(I - \gamma \wh{P}(q_k) \big), \\
        q_{k+1} = q_k -\alpha_k D_k^{-1} \big( q_k - \sbo (q_k)\big),
    \end{array}\right.
    \end{equation*}
    where $\wh{P}(q)$ is the synchronously sampled state-action transition matrix of the Markov chain under the greedy policy~$\pi_q$ w.r.t.~$q$ with elements $[\wh{P}(q)]\big((s,a),(s',a')\big) = 1$ if $s' = \hat{s}^+,\ a' = \pi_q(\hat{s}^+)$ and $=0$ otherwise,
    for each $(s,a),(s',a') \in \set{S}\times\set{A}$, where $\hat{s}^+ \sim \PP(\cdot|s,a)$ is again a \emph{sample} of the next state drawn from the distribution $\PP(\cdot|s,a)$ for the state-action pair $(s,a)$.

%===============================================================================
% \clearpage
\bibliographystyle{apalike} %{siam}
\begin{small}
\bibliography{ref}
\end{small}
%===============================================================================

\end{document}

%% file: AMIN_style.tex
\makeatletter
\def\@settitle{\begin{center}%
		\baselineskip14\p@\relax
		\normalfont\LARGE\bfseries
		\@title
	\end{center}%
}

\def\section{\@startsection{section}{1}%
	\z@{.7\linespacing\@plus\linespacing}{.5\linespacing}%
	{\normalfont\large\bfseries}}

\def\subsection{\@startsection{subsection}{2}%
	\z@{.5\linespacing\@plus.7\linespacing}{.5\linespacing}%
	{\normalfont\bfseries}}

\def\@setauthors{%
  \begingroup
  \def\thanks{\protect\thanks@warning}%
  \trivlist
  \centering\footnotesize \@topsep30\p@\relax
  \advance\@topsep by -\baselineskip
  \item\relax
  \author@andify\authors
  \def\\{\protect\linebreak}%
  \authors%
  \ifx\@empty\contribs
  \else
    ,\penalty-3 \space \@setcontribs
    \@closetoccontribs
  \fi
  \endtrivlist
  \endgroup
}

\makeatother

\usepackage[table, xcdraw, usenames, dvipsnames]{xcolor}
\definecolor{darkblue}{rgb}{0.0, 0.0, 0.45}
\definecolor{darkgreen}{rgb}{0.0, 0.45, 0}
\usepackage[colorlinks	= true,
raiselinks	= true,
linkcolor	= darkblue, %MidnightBlue,
citecolor	= Mahogany,
urlcolor	= darkgreen,
pdfauthor	= {},
pdftitle	= {},
pdfkeywords	= {},
pdfsubject	= {},
plainpages	= false]{hyperref}

\pdfoutput=1
\date{\today}

\usepackage{dsfont,amsfonts,amssymb,amsmath,amsthm}
\usepackage{mathtools, mathrsfs}

\usepackage{graphicx}
\usepackage[font=small,labelfont=rm]{subcaption}
\usepackage[font=small,margin=10pt]{caption}
\usepackage{wrapfig}

\usepackage{multirow}
\usepackage{booktabs}

\usepackage{algorithm,algorithmic}%,algorithmicx}
\usepackage{enumitem}
\usepackage{tikz}
\usepackage{courier} % for embedding the fonts

\usepackage{nicefrac}

\theoremstyle{plain}
\newtheorem{Thm}{Theorem}[section]

\newtheorem{Lem}[Thm]{Lemma}
\newtheorem{Cor}[Thm]{Corollary}

\newtheorem{Rem}[Thm]{Remark}

\DeclareMathOperator*{\argmin}{\arg\!\min}

\newcommand{\R}{\mathbb{R}}

\newcommand{\ra}{\rightarrow}

\newcommand{\wt}{\widetilde}
\newcommand{\wh}{\widehat}
\newcommand{\Let}{\coloneqq}

\newcommand{\tr}{^{\top}}

\newcommand{\norm}[1]{\left\Vert #1 \right\Vert}

\def\ssum{\begingroup\textstyle \sum\endgroup}

\newcommand{\opt}{^\star}
\newcommand{\PP}{\mathds{P}}
\newcommand{\EE}{\mathds{E}}

\newcommand{\set}[1]{\mathcal{#1}}

\DeclareMathOperator{\ord}{\mathcal{O}}